\documentclass[a4paper,10pt]{article}
\usepackage[margin=1in]{geometry}
\usepackage{amssymb, amsmath, amsthm, mathrsfs}
\usepackage{cite}
\usepackage{lineno}
\usepackage{tabls}
\usepackage{multirow}
\usepackage{pifont}
\usepackage{graphicx}

\usepackage{color}
\def\qed{\hfill $\Box$\vspace{0.3cm}}
\def\pf{\noindent{\bf Proof. }}

\newtheorem{lemma}{Lemma}
\newtheorem{theorem}{Theorem}
\newtheorem{remark}{Remark}

\newtheorem{prop}{Proposition}
\newtheorem{definition}{Definition}

\begin{document}

\title{\Large\bf Strong Menger connectedness of augmented $k$-ary $n$-cubes }

\author{
Mei-Mei Gu$^{1,2}$\thanks{This work was partially supported by China Postdoctoral Science Foundation (2018M631322).}
\hspace{.15in}
Jou-Ming Chang$^{3}$\thanks{This work was supported by the grant MOST-107-2221-E-141-001-MY3 from the Ministry of Science and Technology, Taiwan.}
\hspace{.15in}
Rong-Xia Hao$^{1,}$\thanks{This work was partially supported by the National Natural Science Foundation of China (No. 11971054, 11731002) and the 111 Project of China (B16002).
\newline\indent\indent{\footnotesize \emph{E-mail address}: {\tt mmgu@bjtu.edu.cn} (M.-M. Gu), {\tt spade@ntub.edu.tw} (J.-M. Chang),}
{\footnotesize {\tt rxhao@bjtu.edu.cn} (R.-X. Hao)}
}
\\
\\
{\small $^1$ Department of Mathematics, Beijing Jiaotong University, Beijing 100044, P.R. China}\\
{\small  $^2$ Faculty of Mathematics and Physics, Charles University, Prague  11800, Czech Republic}\\
{\small $^3$ Institute of Information and Decision Sciences,}\\
{\small National Taipei University of Business, Taipei 10051, Taiwan}\\
}

\date{}
\maketitle

\begin{abstract}
\baselineskip=11pt
A connected graph $G$ is called strongly Menger (edge) connected if for any two distinct vertices $x,y$ of $G$, there are $\min \{{\rm deg}_G(x), {\rm deg}_G(y)\}$ vertex(edge)-disjoint paths between $x$ and $y$. In this paper, we consider strong Menger (edge) connectedness of the augmented $k$-ary $n$-cube $AQ_{n,k}$, which is a variant of $k$-ary $n$-cube $Q_n^k$. By exploring the topological proprieties of $AQ_{n,k}$, we show that $AQ_{n,3}$ for $n\geq 4$ (resp.\ $AQ_{n,k}$ for $n\geq 2$ and $k\geq 4$) is still strongly Menger connected even when there are $4n-9$ (resp.\ $4n-8$) faulty vertices and $AQ_{n,k}$ is still strongly Menger edge connected even when there are $4n-4$ faulty edges for $n\geq 2$ and $k\geq 3$.
Moreover, under the restricted condition that each vertex has at least two fault-free edges, we show that $AQ_{n,k}$ is still strongly Menger edge connected even when there are $8n-10$ faulty edges for $n\geq 2$ and $k\geq 3$. These results are all optimal in the sense of the maximum number of tolerated vertex (resp.\ edge) faults.


\end{abstract}
\setcounter{page}{1}
\baselineskip=14pt

\section{Introduction}

With continuous advances in technology, a multiprocessor system
may contains hundreds or even thousands of processors that communicate by exchanging messages through an interconnection network. The topology of a network can be represented as a graph.
Among all fundamental properties for interconnection networks, the connectivity and edge connectivity are the major parameters widely discussed for the connection status of networks.

For a connected graph $G$, the {\it connectivity} $\kappa(G)$ is the minimum number of vertices removed to get the graph disconnected or trivial; while the {\it edge connectivity} $\lambda(G)$ is the minimum number of edges removed to get the graph disconnected. Connectivity and edge connectivity are two deterministic measurements for determining the reliability and
fault tolerance of a multiprocessor system.
In contrast with this concept, Menger~\cite{Menger} provided a local point
of view, and defined the connectivity (resp.\ edge connectivity) of any two vertices
as the minimum number of internally vertex-disjoint (resp.\ edge-disjoint) paths between them.

A connected graph $G$ is called {\it strongly Menger (edge) connected} if for any two distinct vertices $x,y$ of $G$, there are $\min\{\text{deg}_G(x),\text{deg}_G(y)\}$ (edge-)disjoint paths between $x$ and $y$.
Parallel routing (i.e., construction of disjoint paths or edge-disjoint paths) has been an important issue in the study of computer networks.
With the continuous increasing in network size, routing in networks with faults has become unavoidable.
Two fault models have been studied for many well-known networks: one is the random fault model, and the other is the conditional fault model (which assumes that the fault distribution is limited). The strong Menger (edge) connectivity of a graph with random faults is defined as follows.

\begin{definition}{\rm
A graph $G$ is called {\it $m$-strongly Menger connected} (resp.\ {\it edge connected}) if $G-F$ remains strongly Menger connected (resp.\ edge connected) for an arbitrary vertex set $F\subseteq V(G)$ (resp.\ edge set $F\subseteq E(G)$) with $|F|\leq m$.
}
\end{definition}

Note that the term $m$-strong Menger connectivity is also referred to as {\it $m$-fault-tolerant maximal local-connectivity} in~\cite{cll15,cct14,sm19}.
Let $\delta(G)$ denote the minimum degree of a graph $G$.
The conditional strong Menger (edge) connectivity of a graph is defined as follows.

\begin{definition}{\rm
A graph $G$ is called {\it $m$-conditional strongly Menger connected} (resp.\ {\it edge connected}) if $G-F$ remains strongly Menger connected (resp.\ edge connected) for an arbitrary vertex set $F\subseteq V(G)$ (resp.\ edge set $F\subseteq E(G)$) with $|F|\leq m$ and $\delta(G-F)\geq 2$.
}
\end{definition}

The study on strong Menger vertex/edge connectedness attracts more and more researchers' attention.
On one hand, the (conditional) strong Menger connectivity of many known networks were explored in literature, for example, star graph~\cite{oc03,lx19t}, hypercubes~\cite{oh04}, folded hypercubes~\cite{yang17}, hypercube-like networks~\cite{lscc}, Cayley graphs generated by transposition trees~\cite{scht09}, augmented cubes~\cite{cct14}, bubble-sort star graph~\cite{cll15}, balanced hypercubes~\cite{lx18} etc.
On the other hand, the (conditional) strong Menger edge connectivity were investigated on hypercubes~\cite{qyang17}, folded hypercubes~\cite{cxu18}, hypercube-like networks~\cite{lx19}, balanced hypercubes~\cite{lx18} etc.
He et al.~\cite{hhc18} and Sabir and Meng~\cite{sm19} presented sufficient conditions for a regular graph to be strongly Menger vertex/edge connected with some restricted conditions.
For more information, please refer to~\cite{hhc18,sm19} and the references therein.

Almost all the known popular classes of networks with strong Menger edge/vertex connectedness are triangle-free.
Recently, augmented $k$-ary $n$-cube is proposed for parallel computing by Xiang and Stewart~\cite{xs11}. An augmented $k$-ary $n$-cube $AQ_{n,k}$ is extended from a $k$-ary $n$-cube $Q_n^k$ in a manner analogous to the extension of an $n$-dimensional hypercube $Q_n$ to an $n$-dimensional augmented $AQ_n$~\cite{cs02} and has many triangles.
Some results about topological properties
and routing problems on augmented $k$-ary $n$-cube
can be found in~\cite{xs11,lz15,gu16,xs09} etc.
In this paper, by exploring and utilizing the structural
properties of $AQ_{n,k}$, we show that $AQ_{n,3}$ (resp.\ $AQ_{n,k}$, $k\geq 4$) is $(4n-9)$-strongly (resp.\ $(4n-8)$-strongly) Menger connected for $n\geq 4$ (resp.\ $n\geq 2$), and $AQ_{n,k}$ is $(4n-4)$-strongly Menger edge connected for $n\geq 2$ and $k\geq 3$.
Moreover, under the restricted condition that each vertex has at least two fault-free edges, we show that $AQ_{n,k}$ is $(8n-10)$-conditional strongly Menger edge connected for $n\geq 2$ and $k\geq 3$.
These results are all optimal with respect to the maximum number of tolerated vertex (resp.\ edge) faults.

\section{Preliminaries}

\subsection{Notations}
Let $G=(V(G),E(G))$ represent an interconnection network,
where a vertex $u\in V(G)$ represents a processor and an edge $(u,v)\in E(G)$
represents a link between vertices $u$ and $v$.

Let $|V(G)|$ be the size of vertex set and $|E(G)|$ be the size of edge set.
Two vertices $u$ and $v$ are {\it adjacent} if $(u,v)\in E(G)$, the vertex $u$ is called a neighbor of $v$, and vice versa. For a vertex $u\in V(G)$, let $N_G(u)$ denote
a set of vertices in $G$ adjacent to $u$, and let $N_{G}[u]=N_G(u)\cup\{u\}$.
The {\it degree} of $u$, denoted by $\text{deg}_G(u)$ (or $d_G(u)$), is the cardinality of $N_{G}(u)$.
For a vertex set $U \subseteq V(G)$, the {\it neighborhood}
of $U$ in $G$ is defined as $N_G(U)=\bigcup\limits_{v\in U}N_{G}(v)-U$.
When no ambiguity arises, we omit the subscript $G$ in the above notations.
For any two vertices $u,v\in V(G)$, we use $cn(G\colon\!u,v)$ to denote the number of common neighbors of $u$ and $v$ in $G$.

A graph $H$ is a {\it subgraph} of a graph $G$ if $V(H)\subseteq V(G)$ and $E(H)\subseteq E(G)$.
The {\em connected components} (simply, {\em component}) of a graph are its maximal connected subgraphs.
For two disjoint vertex sets or subgraphs $H_1$ and $H_2$, we use $E(H_1,H_2)$ to denote the set of edges with one endpoint in $H_1$ and the other in $H_2$.
For a subset $S\subseteq V(G)$ (resp.\ $S\subseteq E(G)$),
we denote $G-S$ the graph obtained from $G$ by removing the vertices (edges) of $S$. In particular, $S$ is called a \emph{vertex cut} (resp.\ {\it edge cut}) of $G$ if $G-S$ is disconnected.
In this case, the biggest component of $G-S$ is called {\it a large component};
the component of $G-S$ which is not the biggest one is called a {\it smaller component}.


Given $x, y\in V(G)$, an {\it $(x,y)$-path}\ of length $k$ is a finite sequence of distinct vertices $v_0,v_1,\ldots, v_k$ such that $x=v_0,y=v_k$, and $(v_i, v_{i+1})\in E(G)$ for $0\leq i\leq k-1$. A set $F\subset V(G)\setminus\{x,y\}$ is an {\it $(x,y)$-cut}\ if
$G-F$ has no $(x, y)$-path. Similarly, a set $F\subseteq E(G)$ is an {\it $(x,y)$-edge cut}\ if $G-F$ has no $(x, y)$-path.

\begin{prop}{\rm(\cite{Menger})}\label{m}
Let $x, y$ be two distinct vertices of a graph $G$.
\begin{enumerate}
\vspace{-0.3cm}
\item [{\rm (1)}] For $(x,y)\notin E(G)$, the minimum size of an $(x,y)$-cut equals the maximum number of disjoint $(x,y)$-paths.
\vspace{-0.2cm}
\item [{\rm (2)}] The minimum size of an $(x, y)$-edge cut equals the maximum number of edge-disjoint $x,y$-paths.
\end{enumerate}
\end{prop}

\subsection{Augmented $k$-ary $n$-cubes}
Let $[n]=\{1,2,\ldots,n\}$ and $[n]_0=\{0,1,\ldots,n-1\}$. Assume that all arithmetics on tuple elements are modulo $k$. Xiang and Stewart~\cite{xs11} gave two equivalent definitions of augmented $k$-ary $n$-cube as follows.

\begin{definition}{\rm(\cite{xs11})}\label{def-aq1}
{\rm Let $n\geq1$ and $k\geq3$ be integers.
The {\it augmented $k$-ary $n$-cube} $AQ_{n,k}$
has $k^n$ vertices, each vertex is labelled by an $n$-bit string
$(a_n,a_{n-1},\ldots,a_2,a_1)$ (or $a_na_{n-1}\cdots a_2a_1$) with $a_i\in[k]_0$ for $i\in[n]$.
There is an edge joining vertex $u=(u_n,u_{n-1},\ldots,u_2,u_1)$
to $v=(v_n,v_{n-1},\ldots,v_2,v_1)$ if and only if one of the following conditions holds.
\begin{enumerate}
\vspace{-0.3cm}
\item [{\rm (1)}] $v_i=u_i-1$ (resp.\ $v_i=u_i+1$)
for some $i\in[n]$ and $u_j=v_j$ for all $j\in[n]\setminus \{i\}$;
and the edge $(u,v)$ is called an {\it $(i,-1)$-edge} (resp.\ {\it $(i,+1)$-edge}).
\vspace{-0.2cm}
\item [{\rm (2)}] for some $2\leq i\leq n$, $v_i=u_i-1$,
$v_{i-1}=u_{i-1}-1$,\ldots, $v_1=u_1-1$ (resp.\ $v_i=u_i+1$,
$v_{i-1}=u_{i-1}+1$,\ldots, $v_1=u_1+1$),
$v_j=u_j$ for all $j> i$; and the edge $(u,v)$ is called an
{\it $(\leq i,-1)$-edge} (resp.\ {\it $(\leq i,+1)$-edge}).
\end{enumerate}
}
\end{definition}

\begin{figure}[htb]
\begin{center}
\includegraphics[width=\textwidth]{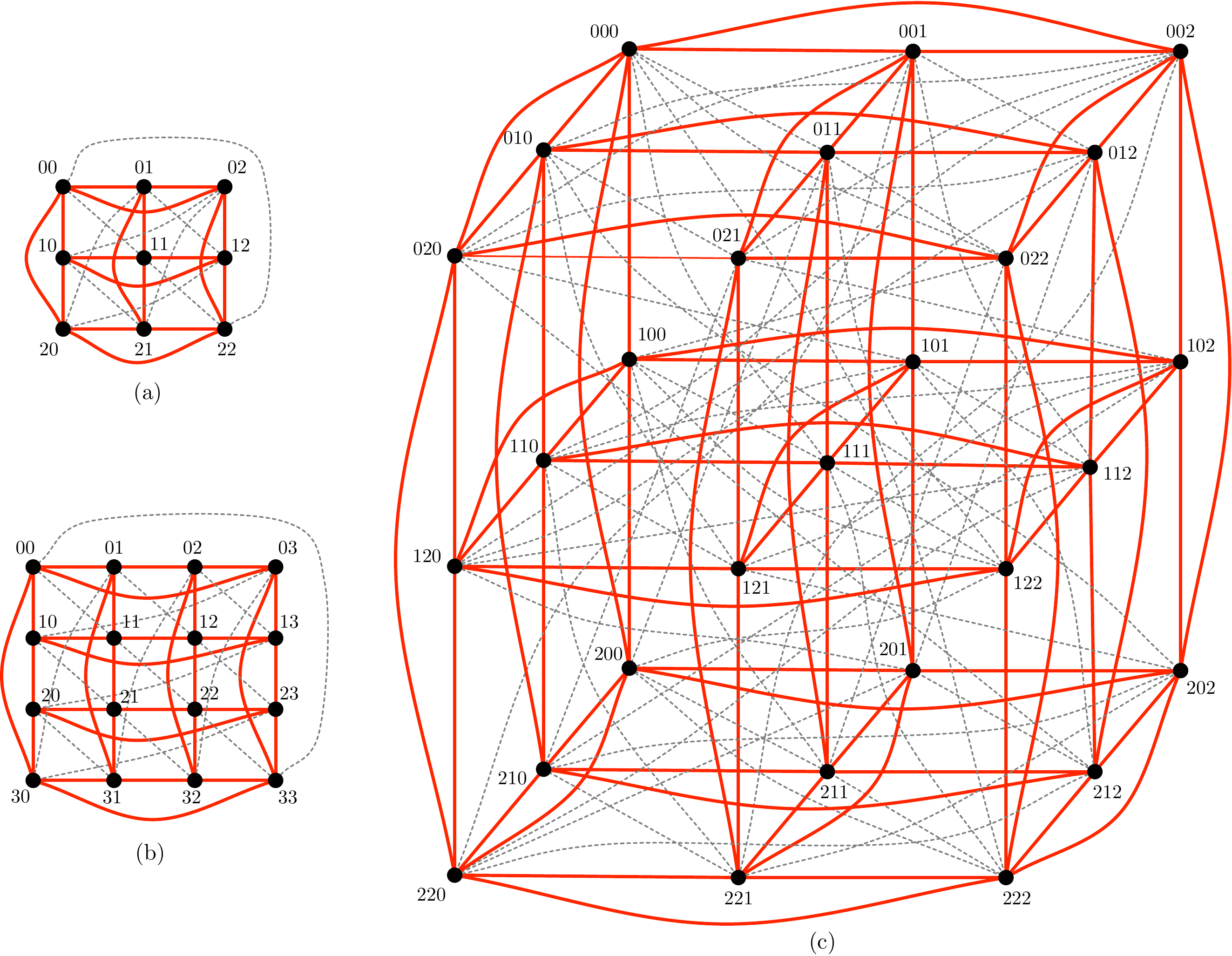}
\caption{(a) $AQ_{2,3}$; (b) $AQ_{2,4}$; (c) $AQ_{3,3}$.}
\label{aq23-24-33}
\end{center}
\end{figure}

In the above definition, edges fulfilled the condition (1) and condition (2) are called {\it traditional edges} and {\it augmented edges}, respectively. Obviously, $AQ_{1,k}$ is a $k$-cycle (i.e., a cycle of length $k$).
Fig.~\ref{aq23-24-33} shows $AQ_{2,3}$, $AQ_{2,4}$ and $AQ_{3,3}$, where bold lines are traditional edges and dashed lines are augmented edges. In fact, the augmented $k$-ary $n$-cube $AQ_{n,k}$ can also be recursively defined as follows.

\begin{definition}{\rm(\cite{xs11})}\label{def-aq2}
{\rm Fix $k\geq 3$, augmented $k$-ary $1$-cube $AQ_{1,k}$ has vertex set
$\{1,2,\ldots,k\}$, and there is an edge joining vertex $u$ to $v$
if and only if $v=u-1$ or $v=u+1$.
Fix $n\geq2$. Take $k$ copies of augmented $k$-ary $(n-1)$-cubes $AQ_{n-1,k}$
and for $i$th copy add an extra number $i$ as the $n$th bit of each vertex
(i.e., all vertices have the same $n$th bit if they are in the same copy of augmented $k$-ary
$(n-1)$-cube). Four more edges are added for each vertex,
namely the $(n,-1)$-edge, $(n,+1)$-edge, $(\leq n,-1)$-edge, $(\leq n,+1)$-edge (as defined in Definition~\ref{def-aq1}).}
\end{definition}

\begin{lemma}{\rm(\cite{xs11})}\label{aq1}
Let $AQ_{n,k}$ be the augmented $k$-ary $n$-cube,
where $n\geq 2$ and $k\geq 3$ are integers.
\begin{enumerate}
\vspace{-0.3cm}
\item [{\rm (1)}] For $i\in[k]_0$,
the subgraph of $AQ_{n,k}$ induced by
vertices with the $n$th bit being $i$, denoted by $AQ_{n,k}^i$,
is a copy of $AQ_{n-1,k}$.
\vspace{-0.2cm}
\item [{\rm (2)}] $AQ_{n,k}$ is vertex-transitive, $AQ_{2,k}$ is edge-transitive
and $\kappa(AQ_{n,k})=\lambda(AQ_{n,k})=4n-2$.
\vspace{-0.2cm}
\item[{\rm (3)}] Every vertex of $AQ_{n,k}^i$ has exactly two neighbors in $AQ_{n,k}^{i+1}$ {\rm(}resp.\ in $AQ_{n,k}^{i-1}${\rm)}, one is connected by a traditional edge and the other is connected by an augmented edge. Thus $|E(AQ_{n,k}^i,AQ_{n,k}^{i+1})|=2k^{n-1}$.
\item[{\rm (4)}] Let $U$ be a subset of $V(AQ_{n,k}^i)$ for $i\in [k]_0$.
Then $|N_{AQ_{n,k}-AQ_{n,k}^i}(U)|\geq 2|U|$.
\end{enumerate}
\end{lemma}

Two adjacent vertices $u\in V(AQ_{n,k}^i)$ and $v\in V(AQ_{n,k}^{i+1})$ are called {\it extra neighbors} to each other. Moreover, $(u,v)$ are called {\it extra edges}. From Lemma~\ref{aq1}(3), every vertex has exactly four distinct extra neighbors.
For convenience, we adopt the following notation to
identify the neighbors of a given vertex in $AQ_{n,k}$.
Let $u=(u_n,u_{n-1},\ldots,u_2,u_1)\in V(AQ_{n,k})$.
For each $2\leq i\leq n$, let

\vspace{-0.2cm}
$u_{(i,-1)}=(u_n,u_{n-1},\ldots,u_i-1,u_{i-1},\ldots,u_2,u_1)$,

\vspace{-0.2cm}
$u_{(i,+1)}=(u_n,u_{n-1},\ldots,u_i+1,u_{i-1},\ldots,u_2,u_1)$,

\vspace{-0.2cm}
$u_{(\leq i,-1)}=(u_n,u_{n-1},\ldots,u_i-1,u_{i-1}-1,\ldots,u_2-1,u_1-1)$,

\vspace{-0.2cm}
$u_{(\leq i,+1)}=(u_n,u_{n-1},\ldots,u_i+1,u_{i-1}+1,\ldots,u_2+1,u_1+1)$.

\begin{remark}{\rm(\cite{lz15})}\label{rmk1}
Let $u$ and $v$ be two distinct vertices in $AQ_{n,k}^i$
such that $u$ and $v$ have common neighbors in $V(AQ_{n,k}^{i+1})\cup V(AQ_{n,k}^{i-1})$.
Then $v=u_{(\leq n-1,-1)}$ or $v=u_{(\leq n-1,+1)}$.
Furthermore, if $v=u_{(\leq n-1,-1)}$, then the common
neighbors of $u$ and $v$ in $V(AQ_{n,k}^{i+1})\cup V(AQ_{n,k}^{i-1})$ are
$u_{(\leq n,-1)}$ and $u_{(n,+1)};$
if $v=u_{(\leq n-1,+1)}$, then the common
neighbors of $u$ and $v$ in $V(AQ_{n,k}^{i+1})\cup V(AQ_{n,k}^{i-1})$
are $u_{(\leq n,+1)}$ and $u_{(n,-1)}$.
\end{remark}

For instance, we consider $u=120$ in $AQ_{3,3}^1$ (see Fig.~\ref{aq23-24-33}(c)). Then, we can check that $u$ and $u_{(\leq 2,-1)}=112$ have common neighbors $u_{(\leq 3,-1)}=012\in V(AQ_{3,3}^0)$ and $u_{(3,+1)}=220\in V(AQ_{3,3}^2)$. Also, $u$ and $u_{(\leq 2,+1)}=101$ have common neighbors $u_{(\leq 3,+1)}=201\in V(AQ_{3,3}^2)$ and $u_{(3,-1)}=020\in V(AQ_{3,3}^0)$.

The following lemma shows
the exact number of common neighbors of any two adjacent vertices in $AQ_{n,k}$.

\begin{lemma}{\rm(\cite{lz15})}\label{aq2}
Let $(u,v)$ be an edge of $AQ_{n,k}$. Then the following assertions hold:
\begin{enumerate}
\vspace{-0.3cm}
\item [{\rm (1)}] $cn(AQ_{2,3}\colon\!u,v)=3$;
\vspace{-0.2cm}
\item [{\rm (2)}] For $k\geq 4$, $cn(AQ_{2,k}\colon\!u,v)=2$;
\vspace{-0.1cm}
\item [{\rm (3)}] For $n\geq 3$,
$cn(AQ_{n,3}\colon\!u,v)=\begin{cases}
3 & \text{if}\ v=u_{(\leq n,-1)}\ \text{or}\ v=u_{(i,-1)}\ \text{for}\ 1\leq i\leq n; \\
5 & \text{if}\ v=u_{(\leq i,-1)}\ \text{for}\ 2\leq i\leq n-1;
\end{cases}$
\vspace{-0.1cm}
\item [{\rm (4)}] For $n\geq 3$ and $k\geq 4$,
$cn(AQ_{n,k}\colon\!u,v)=\begin{cases}
2 & \text{if}\ v=u_{(\leq n,-1)}\ \text{or}\ v=u_{(i,-1)}\ \text{for}\ 1\leq i\leq n; \\
4 & \text{if}\ v=u_{(\leq i,-1)}\ \text{for}\ 2\leq i\leq n-1.\\
\end{cases}$
\vspace{-0.2cm}
\end{enumerate}
The results are similar with assertions {\rm(3)} and {\rm(4)} for $v=u_{(i,+1)}$ and $v=u_{(\leq i,+1)}$.
\end{lemma}


The following lemma shows
the upper bound of the number of common neighbors for any two distinct vertices in $AQ_{n,k}$.

\begin{lemma}{\rm(\cite{lz15})}\label{aq3}
Let $u$ and $v$ be two distinct vertices of $AQ_{n,k}$. Then the following assertions hold:
\begin{enumerate}
\vspace{-0.3cm}
\item [{\rm (1)}] For $n\geq 2$, $cn(AQ_{n,3}\colon\!u,v)\leq 6$;
\vspace{-0.2cm}
\item [{\rm (2)}] For $k\geq 4$, $cn(AQ_{2,k}\colon\!u,v)\leq 2$;
\vspace{-0.2cm}
\item [{\rm (3)}] For $n\geq 3$ and $k\geq 4$, $cn(AQ_{n,k}\colon\!u,v)\leq 4$.
\end{enumerate}
\end{lemma}


%
%

\begin{lemma}{\rm(\cite{gu16})}\label{nes}
Let $AQ_{n,k}$ be the augmented $k$-ary $n$-cubes, where $n\geq 3$ and $k\geq 4$. If
$U$ is a subset of $V(AQ_{n,k})$ with $2\leq |U|\leq 8n-16$, then $|N_{AQ_{n,k}}(U)|\geq 8n-10$.
\end{lemma}

Wang and Zhao~\cite{wz} derived the following result which is useful for our proof.

\begin{lemma}{\rm(\cite{wz})}\label{wz18}
Let $AQ_{n,3}$ be the augmented $3$-ary $n$-cubes, and $F\subseteq V(AQ_{n,3})$ with $|F|\leq 8n-12$. Assume $AQ_{n,3}-F$ is disconnected. Then
\begin{enumerate}
\vspace{-0.3cm}
\item [{\rm (1)}] $AQ_{3,3}-F$ has exactly two
components, one of which is a singleton or a
$3$-cycle, and the vertex set of the $3$-cycle is $\{u,u_{(\leq 2,-1)},u_{(\leq 2,+1)}\}$;
\vspace{-0.2cm}
\item [{\rm (2)}]  for $n\geq 4$, $AQ_{n,3}-F$ has exactly two
components, one of which is a singleton.
\end{enumerate}
\end{lemma}

\section{Main Results}
In this section, we will consider the (conditional) fault-tolerant strong Menger (edge) connectivity of augmented $k$-ary $n$-cubes.
The following result is useful.

\begin{lemma}{\rm(\cite{sm19})}\label{em}
Let $G$ be an $r$-regular, $r$-connected graph with $|V(G)|\geq 2r+1$ and $r\geq 2$. Then $G$ is $f$-strongly Menger connected if, for any $V_f\subseteq V(G)$ with $|V_f|\leq f+r-1$, $G-V_f$ has a component $C$ such that $|V(C)|\geq |V(G)|-|V_f|-1$.
\end{lemma}

\subsection{Strong Menger connectivity of augmented $k$-ary $n$-cubes}

In this section, we consider the strong Menger connectivity of augmented $k$-ary $n$-cubes $AQ_{n,k}$. We will prove that $AQ_{n,3}$ is $(4n-9)$-strongly Menger connected but not $(4n-8)$-strongly Menger connected for $n\geq 3$, and $AQ_{n,k}$ is $(4n-8)$-strongly Menger connected but not $(4n-7)$-strongly Menger connected for $n\geq 3$ and $k\geq 4$.

\begin{lemma}\label{vertex}
Let $F$ be an arbitrary set of vertices in $AQ_{n,k}$.
\begin{enumerate}
\vspace{-0.3cm}
\item [{\rm (1)}] If $|F|\leq 8n-12$ for $n\geq 4$, then $AQ_{n,3}-F$ has a component $C$ such that $|V(C)|\geq |V(AQ_{n,3})|-|F|-1$.
\vspace{-0.2cm}
\item [{\rm (2)}] If $|F|\leq 8n-11$ for $n\geq 2$ and $k\geq 4$, then $AQ_{n,k}-F$ has a component $C$ such that $|V(C)|\geq |V(AQ_{n,k})|-|F|-1$.
\end{enumerate}
\end{lemma}

\pf Let $C$ be the large component of $AQ_{n,k}-F$.
First we consider $k=3$, $n\geq 4$ and $|F|\leq 8n-12$. By Lemma~\ref{wz18}(2), $AQ_{n,3}-F$ is connected or has exactly two components, one of which is a singleton. Clearly, if $AQ_{n,k}-F$ is connected, then $V(C)=V(AQ_{n,k}-F)$ and $|V(C)|=|V(AQ_{n,3})|-|F|$. If $AQ_{n,3}-F$ is disconnected, then  $|V(C)|=|V(AQ_{n,3})|-|F|-1$.

Now assume that $n\geq 2$, $k\geq 4$ and $|F|\leq 8n-11$. The proof is by induction on $n$. If $n=2$, then $|F|\leq 5<\kappa(AQ_{2,k})=6$. Thus, $AQ_{2,k}-F$ is connected. It leads to $V(C)=V(AQ_{2,k}-F)$ and $|V(C)|=|V(AQ_{2,k})|-|F|$.
In what follows, we assume that $n\geq 3$ and the result holds for $AQ_{n-1,k}$.
Recall that $AQ_{n,k}$ contains $k$ disjoint copies of $AQ_{n-1,k}$, say $AQ_{n,k}^i$, $i\in [k]_0$. Let $F_i=F\cap V(AQ_{n,k}^i)$ and $f_i=|F_i|$ for $i\in [k]_0$.
Let $I=\{i\in [k]_0\colon\,AQ_{n,k}^i-F_i\ \text{is disconnected}\}$ and $J=[k]_0\setminus I$.
In addition, we adopt the following notations:
\[
F_I=\bigcup_{i\in I}F_i,\
F_J=\bigcup_{j\in J}F_j,\ \text{and}\
AQ_{n,k}^J=\bigcup_{j\in J}AQ_{n,k}^j.
\]
By Lemma~\ref{aq1}(2), $f_i\geq 4(n-1)-2=4n-6$ for $i\in I$.
Since $|F|\leq 8n-11$, we have $|I|\leq 2$.
We consider the following cases.

Case 1: $|I|=0$.

For all $j\in[k]_0$, $AQ_{n,k}^j-F_j$ is connected. By Lemma~\ref{aq1}(3), there are $2k^{n-1}$ edges between subgraphs $AQ_{n,k}^j$ and $AQ_{n,k}^{j+1}$. Since $|F|\leq 8n-11<2k^{n-1}$
for $n\geq 3$ and $k\geq 4$, there is a fault-free edge between $AQ_{n,k}^j-F_j$ and $AQ_{n,k}^{j+1}-F_{j+1}$ for each $j\in[k]_0$, it implies that $AQ_{n,k}-F$ is connected. Let $C=AQ_{n,k}-F$. Clearly, $|V(C)|=|V(AQ_{n,k})|-|F|$.

Case 2: $|I|=1$.

Without loss of generality, assume that $I=\{0\}$. By Lemma~\ref{aq1}(2), $f_0\geq 4n-6$. For $j\in [k]_0\setminus \{0\}$, $AQ_{n,k}^j-S_j$ is connected.
By Lemma~\ref{aq1}(3), there are $2k^{n-1}$ edges between subgraphs $AQ_{n,k}^j$ and $AQ_{n,k}^{j+1}$. For each $j,j+1\in[k]_0\setminus I$, since $2(f_j+f_{j+1})\leq 2|F|\leq 2(8n-11)<2k^{n-1}$
for $n\geq 3$ and $k\geq 4$, there is a fault-free edge between $AQ_{n,k}^j-F_j$ and $AQ_{n,k}^{j+1}-F_{j+1}$. It leads to $AQ_{n,k}^J-F_J$ is connected.

Case 2.1: $4n-6\leq f_0\leq 8n-19$.

Since $f_0\leq 8n-19=8(n-1)-11$, by induction hypothesis on $AQ_{n,k}^0$, there exists a component, say $H_0$ in $AQ_{n,k}^0$, such that $|V(H_0)|\geq |V(AQ_{n,k}^0)|-f_0-1=k^{n-1}-f_0-1$. Since
there are $2\times2k^{n-1}$ edges between subgraphs $AQ_{n,k}^0$ and $AQ_{n,k}^J$ and
$2(f_0+f_1+1)+2(f_0+f_{k-1}+1)\leq 4|F|+4\leq 4(8n-11)+4<2\times2k^{n-1}$
for $n\geq 3$ and $k\geq 4$, $H_0$ is connected to $AQ_{n,k}^J-F_J$. Let $C$ be the component induced by the vertex set $V(H_0)\cup V(AQ_{n,k}^J-F_J)$. Then $|V(C)|\geq |V(AQ_{n,k})|-|F|-1$.

Case 2.2: $8n-18\leq f_0\leq 8n-11$.

In this case, $|F_J|=|F|-f_0\leq 7$. Let $H_0$ be the large component of $AQ_{n,k}^0-F_0$ and $M=AQ_{n,k}^0-F_0-V(H_0)$. Obviously, $N_{AQ_{n,k}-AQ_{n,k}^0}(V(M))\subseteq F_J$ and $N_{AQ_{n,k}^0}V(M)\subseteq F_0$.
By Lemma~\ref{aq1}(4), $|N_{AQ_{n,k}-AQ_{n,k}^0}V(M)|\geq 2|V(M)|$.
It leads to $2|V(M)|\leq |F_J|\leq 7$, so $|V(M)|\leq 3$.
If $2\leq |V(M)|\leq 3$, by Lemma~\ref{nes}, $|N_{AQ_{n,k}}(V(M))|\geq 8n-10$.
It leads to
\[
8n-11\geq |F|=f_0+|F_J|\geq |N_{AQ_{n,k}^0}V(M)|+|N_{AQ_{n,k}-AQ_{n,k}^0}(V(M))|=|N_{AQ_{n,k}}(V(M))|\geq 8n-10,
\]
a contradiction. Thus, $|V(M)|=1$ and $|V(H_0)|=|V(AQ_{n,k}^0)|-f_0-1$. One can see that $H_0$ is connected to $AQ_{n,k}^J-F_J$ by the similar argument as Case 2.1. Let $C$ be the component induced by the vertex set $V(H_0)\cup V(AQ_{n,k}^J-F_J)$. Then $|V(C)|\geq |V(AQ_{n,k})|-|F|-1$.

Case 3: $|I|=2$.

Without loss of generality, assume that $I=\{0,t\}$, where $1\leq t\leq k-1$. By Lemma~\ref{aq1}(2), $f_0,  f_t\geq 4n-6$. Since $|F|\leq 8n-11$, we have $|F_J|\leq 1$. For $j\in [k]_0\setminus \{0,t\}$, $AQ_{n,k}^j-F_j$ is connected.
We consider the following cases.

Case 3.1: $t=1$ or $k-1$.

Without loss of generality, assume that $t=1$.
Note that $|F_J|\leq 1$ and $2k^{n-1}>2|F_J|$ for $n\geq 3$.
For $j,j+1\in [k]_0\setminus \{0,1\}$, there is a fault-free edge between $AQ_{n,k}^j-F_j$ and $AQ_{n,k}^{j+1}-F_{j+1}$. It leads to $AQ_{n,k}^J-F_J$ is connected.
Since every vertex of $AQ_{n,k}^0$ (resp.\ $AQ_{n,k}^1$) has two extra
neighbors in $AQ_{n,k}^J$ and $|F_J|\leq 1$,
any component of $AQ_{n,k}^0-F_0$ (resp.\ $AQ_{n,k}^1-F_1$) is connected to $AQ_{n,k}^J-F_J$. Let $C=AQ_{n,k}-F$. Clearly, $|V(C)|=|V(AQ_{n,k})|-|F|$.

Case 3.2: $2\leq t\leq k-2$.

For $1\leq m\leq t-1$ or $t+1\leq m\leq k-1$, $AQ_{n,k}^m-F_m$ is connected.
By the similar argument as Case~3.1, those $(AQ_{n,k}^m-F_m)$'s for $1\leq m\leq t-1$ and $t+1\leq m\leq k-1$ belong to the components, say $C_1$ and $C_2$, respectively, of $AQ_{n,k}-F$.
Since every vertex of $AQ_{n,k}^0$ (resp.\ $AQ_{n,k}^t$) has four extra
neighbors in $AQ_{n,k}-AQ_{n,k}^0-AQ_{n,k}^t$
and $|F_J|\leq 1$, any component of $AQ_{n,k}^0-F_0$
is connected to both $AQ_{n,k}^1-F_1$ (which is part of $C_1$) and $AQ_{n,k}^{k-1}-F_{k-1}$ (which is part of $C_2$). This implies that $C_1=C_2$ (i.e., $C_1$ and $C_2$ are the same component).
By a similar discussion, any component of $AQ_{n,k}^t-F_t$ is contained in both $C_1$ and $C_2$.
This implies that $AQ_{n,k}-F$ has a large component, say $C$, and $C=C_1=C_2$.
It leads to $|V(C)|=|V(AQ_{n,k})|-|F|$.
\qed

Since $AQ_{n,k}$ is $(4n-2)$-regular and $(4n-2)$-connected by Lemma~\ref{aq1}(2), combining Lemma~\ref{em} and Lemma~\ref{vertex}, we have the following result.

\begin{theorem}\label{th1}
Let $AQ_{n,k}$ be the augmented $k$-ary $n$-cube. Then
\begin{enumerate}
\vspace{-0.3cm}
\item [{\rm (1)}] $AQ_{n,3}$ is $(4n-9)$-strongly Menger edge connected for $n\geq 4$.
\vspace{-0.2cm}
\item [{\rm (2)}] $AQ_{n,k}$ is $(4n-8)$-strongly Menger edge connected for $n\geq 2$ and $k\geq 4$.
\end{enumerate}
\end{theorem}

\begin{figure}[htb]
\begin{center}
\includegraphics[width=3.1in]{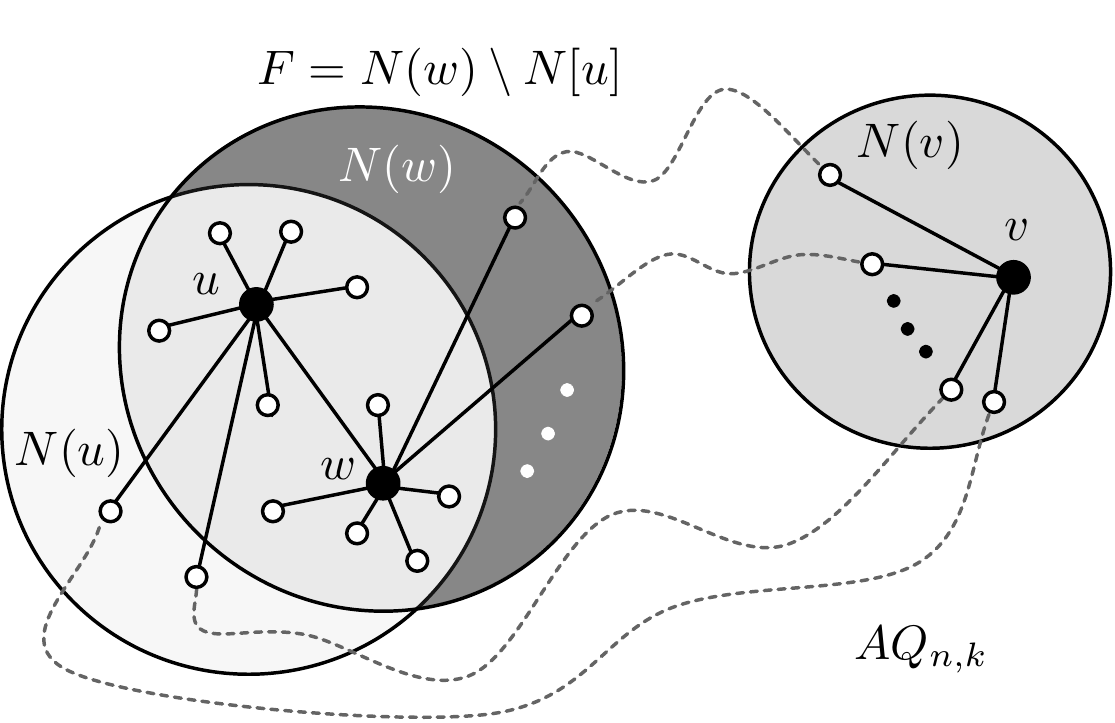}
\caption{Illustration for Remark~\ref{rmk2}}
\label{remark2}
\end{center}
\end{figure}

\begin{remark}\label{rmk2}
{\rm
$AQ_{n,3}$ is not $(4n-8)$-strongly Menger connected for $n\geq 3$ and $AQ_{n,k}$ is not $(4n-7)$-strongly Menger connected for $n\geq 3$ and $k\geq 4$. See Fig.~\ref{remark2} for an illustration. Let $(u,w)\in E(AQ_{n,k})$ such that $u=w_{(\leq i,-1)}$ for some $1\leq i\leq n-1$ and let $F=N(w)\setminus N[u]$ be a faulty subset of vertices in $AQ_{n,k}$ (i.e., vertices in the darkest area of Fig.~\ref{remark2}). By Lemma~\ref{aq2}(3), $|F|=(4n-2)-5-1=4n-8$ if $n\geq 3$ and $k=3$, and by Lemma~\ref{aq2}(4), $|F|=(4n-2)-4-1=4n-7$ if $n\geq 3$ and $k\geq 4$. We now consider a vertex $v\in V(AQ_{n,k})\setminus (N[N(w)\cup N(u)])$. Obviously, ${\rm deg}_{AQ_{n,k}-F}(u) ={\rm deg}_{AQ_{n,k}-F}(v)=4n-2$. Since $w\in N(u)$ and some neighbors of $w$ are in $F$, the vertices $u$ and $v$ are not connected with $4n-2$ vertex-disjoint paths in $AQ_{n,k}-F$. Thus, the results of Theorem~\ref{th1} are optimal in the sense that the number of faulty vertices cannot be increased.
}
\end{remark}

\subsection{Strong Menger edge connectivity of augmented $k$-ary $n$-cubes}

In this section, we consider the (conditional) strongly Menger edge connectivity of augmented $k$-ary $n$-cubes.

In the following, let $S$ be an arbitrary set of edges in $AQ_{n,k}$.
Note that $AQ_{n,k}$ contains $k$ disjoint copies of $AQ_{n-1,k}$, say $AQ_{n,k}^i$, $i\in [k]_0$. Let $S_i=S\cap E(AQ_{n,k}^i)$ and $s_i=|S_i|$ for $i\in [k]_0$.
Let $I=\{i\in [k]_0\colon\,AQ_{n,k}^i-S_i\ \text{is disconnected}\}$ and $J=[k]_0\setminus I$.
In addition, we adopt the following notations:
\[
S_I=\bigcup_{i\in I}S_i,\
S_J=\bigcup_{j\in J}S_j,\
AQ_{n,k}^J=\bigcup_{j\in J}AQ_{n,k}^j,\ \text{and}\
s_c=|S|-\sum_{i\in [k]_0}s_i.
\]

First, we provide two useful lemmas as follows.

\begin{lemma}\label{edge1}
Let $S$ be an arbitrary set of edges in $AQ_{n,k}$ for $n\geq 2$ and $k\geq 3$. If $|S|\leq 8n-7$, then there exists a component $H$ in $AQ_{n,k}-S$ such that $|V(H)|\geq |V(AQ_{n,k})|-1$.
\end{lemma}


\pf Let $H$ be the large component of $AQ_{n,k}-S$. The proof is by induction on $n$. For $n=2$, the proof is shown in the Appendix A. In what follows, we assume that $n\geq 3$ and $k\geq 3$ and the result holds for $AQ_{n-1,k}$.
Recall that $I=\{i\in [k]_0\colon\,AQ_{n,k}^i-S_i\ \text{is disconnected}\}$ and $J=[k]_0\setminus I$.
By Lemma~\ref{aq1}(2), $s_i\geq 4(n-1)-2=4n-6$ for $i\in I$.
Since $|S|\leq 8n-7$, we have $|I|\leq 2$ when $n\geq 3$. The following cases should be considered.

Case 1: $|I|=0$.

For all $j\in[k]_0$, $AQ_{n,k}^j-S_j$ is connected. By Lemma~\ref{aq1}(3), there are $2k^{n-1}$ edges between subgraphs $AQ_{n,k}^j$ and $AQ_{n,k}^{j+1}$. Since $|S|\leq 8n-7<2k^{n-1}$ for $n\geq 3$ and $k\geq 3$,  there is a fault-free edge between $AQ_{n,k}^j-S_j$ and $AQ_{n,k}^{j+1}-S_{j+1}$ for each $j\in[k]_0$, it implies that $AQ_{n,k}-S$ is connected. Thus, $|V(H)|=|V(AQ_{n,k})|$.

Case 2: $|I|=1$.

Without loss of generality, assume that $I=\{0\}$. By Lemma~\ref{aq1}(2), $s_0\geq 4(n-1)-2=4n-6$. For $j\in [k]_0\setminus \{0\}$, $AQ_{n,k}^j-S_j$ is connected.
By Lemma~\ref{aq1}(3), there are $2k^{n-1}$ edges between subgraphs $AQ_{n,k}^j$ and $AQ_{n,k}^{j+1}$. Since $s_c\leq|S|-s_0\leq(8n-7)-(4n-6)=4n-1<2k^{n-1}$ for $n\geq 3$ and $k\geq 3$, there is a fault-free edge between $AQ_{n,k}^j-S_j$ and $AQ_{n,k}^{j+1}-S_{j+1}$ for each $j,j+1\in J$. It leads to $AQ_{n,k}^J-S_J$ is connected.

Case 2.1: $4n-6\leq s_0\leq 8n-15$.

Since $s_0\leq 8n-15=8(n-1)-7$, by induction hypothesis on $AQ_{n,k}^0$, there exists a component, say $H_0$ in $AQ_{n,k}^0$, such that $|V(H_0)|\geq |V(AQ_{n,k}^0)|-1=k^{n-1}-1$. Since every vertex of $H_0$ has exactly four distinct extra neighbors, we have
$s_c\leq|S|-s_0\leq 4n-1<4\times(k^{n-1}-1)$
for $n\geq 3$ and $k\geq 3$, and thus $H_0$ is connected to $AQ_{n,k}^J-S_J$. Let $H$ be the component induced by the vertex set $V(H_0)\cup V(AQ_{n,k}^J-S_J)$. Then $|V(H)|\geq |V(AQ_{n,k})|-1$, as desired.

Case 2.2:  $8n-14\leq s_0\leq 8n-7$.

In this case, we have $s_c\leq |S|-s_0\leq (8n-7)-(8n-14)=7$. Since every vertex
in $AQ_{n,k}^0$ has exactly four distinct extra neighbors, at most one vertex in
$AQ_{n,k}^0$ are not connected with $AQ_{n,k}^J-S_J$. This shows that $|V(H)|\geq |V(AQ_{n,k})|-1$, as desired.

Case 3:  $|I|=2$.

We consider the following two cases according to $k=3$ or not.

Case 3.1: $k=3$.

Without loss of generality, assume that $I=\{0,1\}$ and $s_0\leq s_1$. Then, $AQ_{n,3}^J-S_J=AQ_{n,3}^2-S_2$ is connected. Since $|S|\leq 8n-7$, if $s_0\geq 8n-14$, then $s_c\leq|S|-s_0-s_1\leq (8n-7)-2(8n-14)=21-8n$, which contradicts that $n\geq 3$. Indeed, by Lemma~\ref{aq1}(2), we have $4n-6\leq s_0\leq s_1\leq |S|-(4n-6)\leq 4n-1$ and $s_c\leq|S|-s_0-s_1\leq (8n-7)-2(4n-6)=5$. Thus, we consider the following two situations.

Case 3.1.1: $4n-6\leq s_0\leq s_1\leq 8n-15=8(n-1)-7$.

Note that $AQ_{n,3}^0-S_0$ (resp.\ $AQ_{n,3}^1-S_1$) is disconnected. By induction hypothesis on $AQ_{n,3}^0$ (resp.\ $AQ_{n,3}^1$), $AQ_{n,3}^0-S_0$ (resp.\ $AQ_{n,3}^1-S_1$) has a large component and a singleton, say $x_0$ (resp.\ $x_1$). Note that each singleton has exactly two distinct extra neighbors in $AQ_{n,3}^2$. Since $s_c\leq5<2k^{n-1}-2$ for $n\geq 3$ and $k=3$, $AQ_{n,3}^0-S_0-\{x_0\}$ (resp.\ $AQ_{n,3}^1-S_1-\{x_1\}$) is connected to $AQ_{n,3}^2-S_2$.
Let $M$ be the union of smaller components of $AQ_{n,3}-S$. Clearly, $V(M)\subseteq \{x_0,x_1\}$.
If $|V(M)|=2$, then $s_c\geq |N_{AQ_{n,3}^1\cup AQ_{n,3}^2}(x_0)|+|N_{AQ_{n,3}^2}(x_1)|\geq 6$, a contradiction occurs. Thus, $|V(M)|\leq 1$.
This implies that $AQ_{n,3}-S$ has a large component $H$ and smaller components which contain at most one  vertices in total. It leads to $|V(H)|\geq |V(AQ_{n,3})|-1$.

Case 3.1.2: $4n-6\leq s_0\leq 8n-15$ and $8n-14\leq s_1\leq 4n-1$.

In this case, we have $(4n-6)+(8n-14)\leq s_0+s_1\leq |S|\leq 8n-7$. It implies that $n=3$.
Thus, we have $|S|\leq 17$, $6\leq s_0\leq 9$, $10\leq s_1\leq 11$ and $s_c\leq|S|-s_1-s_2\leq 1$. Since every vertex in $AQ_{n,3}^0$ (resp.\ $AQ_{n,3}^1$) has two distinct extra neighbors in $AQ_{n,3}^2$, any component in $AQ_{n,3}^0-S_0$ (resp.\ $AQ_{n,3}^1-S_1$) is connected to $AQ_{n,3}^2-S_2$. This implies  $AQ_{n,3}-S$ is connected and  $|V(H)|=|V(AQ_{n,3})|$.

Case 3.2: $k\geq 4$.

Without loss of generality, assume that $0\in I$. We consider the following two cases according to the value of $I\setminus\{0\}$.

Case 3.2.1: $I=\{0,1\}$ or $I=\{0,k-1\}$.

Without loss of generality, assume that $I=\{0,1\}$ and $s_0\leq s_1$. For $j\in [k]_0\setminus\{0,1\}$, $AQ_{n,k}^j-S_j$ is connected. By Lemma~\ref{aq1}(3), since $s_c\leq|S|-s_0-s_1\leq (8n-7)-2(4n-6)=5<2k^{n-1}$ for $n\geq 3$ and $k\geq 4$, there exists a fault-free edge joining $AQ_{n,k}^j-S_j$ and $AQ_{n,k}^{j+1}-S_{j+1}$ for $j,j+1\in [k]_0\setminus\{0,1,k-1\}$. It leads to $AQ_{n,k}^J-S_J$ is connected. By Lemma~\ref{aq1}(2), since $|S|\leq 8n-7$, we have $4n-6\leq s_0\leq s_1\leq |S|-(4n-6)\leq 4n-1$. The following two cases should be considered.

Case 3.2.1a: $4n-6\leq s_0\leq s_1\leq 8n-15$.

Case 3.2.1b: $4n-6\leq s_0\leq 8n-15$ and $8n-14\leq s_1\leq 4n-1$.

Since the discussions for Case 3.2.1a and Case 3.2.1b are similar to those of Case 3.1.1 and Case 3.1.2 respectively, the details are omitted.

Case 3.2.2: $I=\{0,t\}$, where $2\leq t\leq k-2$.

For $1\leq m\leq t-1$ or $t+1\leq m\leq k-1$, $AQ_{n,k}^m-S_m$ is connected.
By the similar argument as Case~3.2.1, those $(AQ_{n,k}^m-S_m)$'s for $1\leq m\leq t-1$ and $t+1\leq m\leq k-1$ belong to the components, say $C_1$ and $C_2$, respectively, of $AQ_{n,k}-S$.
Without loss of generality, assume that $s_0\leq s_t$. We consider the following cases.

Case 3.2.2a: $4n-6\leq s_0\leq s_t\leq 8n-15=8(n-1)-7$.

In this case, $AQ_{n,k}^0-S_0$ (resp.\ $AQ_{n,k}^t-S_t$) is disconnected. By induction hypothesis on $AQ_{n,k}^0$ (resp.\ $AQ_{n,k}^t$),
$AQ_{n,k}^0-S_0$ (resp.\ $AQ_{n,k}^t-S_t$) has exactly two components: a large component, say $H_0$ (resp.\ $H_t$), and a singleton. Note that every vertex in $AQ_{n,k}^0$ has exactly two distinct extra neighbors in $AQ_{n,k}^1$ (resp.\ $AQ_{n,k}^{k-1}$).
Since $s_c\leq|S|-s_0-s_t\leq 5<2k^{n-1}-4$ for $n\geq 3$ and $k\geq 4$,
there is at least a fault-free edge between $H_0$ and $AQ_{n,k}^1-S_1$ (resp.\ $AQ_{n,k}^{k-1}-S_{k-1}$). This implies that $H_0$ is connected to $AQ_{n,k}^1-S_1$ (which is part of $C_1$) and is connected to $AQ_{n,k}^{k-1}-S_{k-1}$ (which is part of $C_2$). It follows that $H_0$ is contained in both $C_1$ and $C_2$. By a similar discussion, $H_t$ is contained in both $C_1$ and $C_2$.

Let $M$ be the union of smaller components of $AQ_{n,k}-S$. Clearly, $V(M)\subseteq \{x_0,x_t\}$.
If $|V(M)|=2$, then $s_c\geq |N_{AQ_{n,k}-AQ_{n,k}^0}(x_0)|+|N_{AQ_{n,k}-AQ_{n,k}^t}(x_t)|\geq 8$, a contradiction occurs. Thus, $|V(M)|\leq 1$.
This implies that $AQ_{n,k}-S$ has a large component $H$ and smaller components which contain at most one  vertices in total. So $|V(H)|\geq |V(AQ_{n,k})|-1$.

Case 3.2.2b: $4n-6\leq s_0\leq 8n-15$ and $8n-14\leq s_t\leq 4n-1$.

In this case, we have $(4n-6)+(8n-14)\leq s_0+s_t\leq |S|\leq 8n-7$. It implies that $n=3$.
Thus, $|S|\leq 17$, $6\leq s_0\leq 9$, $10\leq s_t\leq 11$ and $s_c\leq|S|-s_1-s_t\leq 1$.  Since every vertex in $AQ_{n,k}^0$ (resp.\ $AQ_{n,k}^t$) has four distinct extra neighbors in $AQ_{n,k}-AQ_{n,k}^0-AQ_{n,k}^t$, any component in $AQ_{n,k}^0-S_0$ (resp.\ $AQ_{n,k}^t-S_t$) is connected to $AQ_{n,k}^1-S_1$ (resp.\ $AQ_{n,k}^{t-1}-S_{t-1}$) (which is part of $C_1$) and is connected to $AQ_{n,k}^{k-1}-S_{k-1}$ (resp.\ $AQ_{n,k}^{t+1}-S_{t+1}$) (which is part of $C_2$).
This implies $AQ_{n,k}-S$ is connected and thus $|V(H)|=|V(AQ_{n,k})|$.
\qed

\begin{lemma}\label{edge2}
Let $S$ be an arbitrary set of edges in $AQ_{n,k}$ for $n\geq 2$ and $k\geq 3$. If $|S|\leq 12n-13$, then there exists a component $H$ in $AQ_{n,k}-S$ such that $|V(H)|\geq |V(AQ_{n,k})|-2$.
\end{lemma}
\pf Let $H$ be the large component of $AQ_{n,k}-S$. The proof is by induction on $n$. For $n=2$, the proof is shown in the Appendix B.
In what follows, we assume that $n\geq 3$ and $k\geq 3$ and the result holds for $AQ_{n-1,k}$.
Recall that $I=\{i\in [k]_0\colon\,AQ_{n,k}^i-S_i\ \text{is disconnected}\}$ and $J=[k]_0\setminus I$.
By Lemma~\ref{aq1}(2), $s_i\geq 4(n-1)-2=4n-6$ for $i\in I$.
Since $|S|\leq 12n-13$, we have $|I|\leq 3$ when $n\geq 3$.
We consider the following cases.

Case 1: $|I|=0$.

For all $j\in[k]_0$, $AQ_{n,k}^j-S_j$ is connected. By Lemma~\ref{aq1}(3), there are $2k^{n-1}$ edges between subgraphs $AQ_{n,k}^j$ and $AQ_{n,k}^{j+1}$. Since $|S|\leq 12n-13<2\times(2k^{n-1})$ for $n\geq 3$ and $k\geq 3$, there exists at most one integer, say $i\in[k]_0$, such that all the edges between $AQ_{n,k}^i$ and $AQ_{n,k}^{i+1}$ are faulty. Since there is a fault-free edge between $AQ_{n,k}^j-S_j$ and $AQ_{n,k}^{j+1}-S_{j+1}$ for each $j\in[k]_0\setminus\{i\}$, it implies that $AQ_{n,k}-S$ is connected. Thus, $|V(H)|=|V(AQ_{n,k})|$.

Case 2: $|I|=1$.

Without loss of generality, assume that $I=\{0\}$. By Lemma~\ref{aq1}(2), $s_0\geq 4(n-1)-2=4n-6$. For $j\in [k]_0\setminus \{0\}$, $AQ_{n,k}^j-S_j$ is connected.
By Lemma~\ref{aq1}(3), there are $2k^{n-1}$ edges between subgraphs $AQ_{n,k}^j$ and $AQ_{n,k}^{j+1}$. Since $s_c\leq|S|-s_0\leq(12n-13)-(4n-6)=8n-7<2k^{n-1}$ for $n\geq 3$ and $k\geq 3$, there is a fault-free edge between $AQ_{n,k}^j-S_j$ and $AQ_{n,k}^{j+1}-S_{j+1}$ for each $j,j+1\in J$. It leads to $AQ_{n,k}^J-S_J$ is connected.

Case 2.1: $4n-6\leq s_0\leq 12n-25$.

Since $s_0\leq 12n-25=12(n-1)-13$, by induction hypothesis on $AQ_{n,k}^0$, there exists a component, say $H_0$ in $AQ_{n,k}^0$, such that $|V(H_0)|\geq |V(AQ_{n,k}^0)|-2=k^{n-1}-2$. Since every vertex of $H_0$ has exactly four distinct extra neighbors, we have $s_c\leq|S|-s_0\leq 8n-7<4\times(k^{n-1}-2)$ for $n\geq 3$ and $k\geq 3$, and thus $H_0$ is connected to $AQ_{n,k}^J-S_J$. Let $H$ be the component induced by the vertex set $V(H_0)\cup V(AQ_{n,k}^J-S_J)$. Then $|V(H)|\geq |V(AQ_{n,k})|-2$.

Case 2.2:  $12n-24\leq s_0\leq 12n-13$.

In this case, we have $s_c\leq |S|-s_0\leq (12n-13)-(12n-24)=9$. Since every vertex
in $AQ_{n,k}^0$ has exactly four distinct extra neighbors, at most two vertices in
$AQ_{n,k}^0$ are not connected with $AQ_{n,k}^J-S_J$. This shows that $|V(H)|\geq |V(AQ_{n,k})|-2$, as desired.

Case 3: $|I|=2$.

We consider the following two cases according to $k=3$ or not.

Case 3.1: $k=3$.

Without loss of generality, assume that $I=\{0,1\}$ and $s_0\leq s_1$. Then, $AQ_{n,3}^J-S_J=AQ_{n,3}^2-S_2$ is connected. By Lemma~\ref{aq1}(2), since $|S|\leq 12n-13$, we have $4n-6\leq s_0\leq s_1\leq |S|-(4n-6)\leq 8n-7$.

Case 3.1.1: $4n-6\leq s_0\leq s_1\leq 8n-15=8(n-1)-7$.

Note that $AQ_{n,3}^0-S_0$ (resp.\ $AQ_{n,3}^1-S_1$) is disconnected. By Lemma~\ref{edge1}, $AQ_{n,3}^0-S_0$ (resp.\ $AQ_{n,3}^1-S_1$) has a large component and a singleton, say $x_0$ (resp.\ $x_1$). Note that the singleton has exactly two distinct extra neighbors in $AQ_{n,3}^2$. Since $s_c\leq|S|-s_0-s_1\leq (12n-13)-2(4n-6)=4n-1<2k^{n-1}-2$ for $n\geq 3$ and $k=3$, $AQ_{n,3}^0-S_0-\{x_0\}$ (resp.\ $AQ_{n,3}^1-S_1-\{x_1\}$) is connected to $AQ_{n,3}^2-S_2$. This implies that $AQ_{n,3}-S$ has a large component $H$ and smaller components which contain at most two vertices in total. It leads to $|V(H)|\geq |V(AQ_{n,k})|-2$.

Case 3.1.2: $4n-6\leq s_0\leq 8n-15$ and $8n-14\leq s_1\leq 8n-7$.

In this case, we have $s_c\leq|S|-s_1-s_2\leq (12n-13)-(4n-6)-(8n-14)=7$. Since every vertex in $AQ_{n,3}^0$ (resp.\ $AQ_{n,3}^1$) has four distinct extra neighbors, any component with more than two vertices in $AQ_{n,3}^0-S_0$ (resp.\ $AQ_{n,3}^1-S_1$) is connected to $AQ_{n,3}^2-S_2$. This implies that only a component with a singleton in $AQ_{n,3}^0-S_0$ (resp.\ $AQ_{n,3}^1-S_1$) can be disconnected with $AQ_{n,3}^2-S_2$. Thus, $AQ_{n,3}-S$ has a large component $H$ and smaller components which contain at most two vertices in total. It leads to $|V(H)|\geq |V(AQ_{n,3})|-2$.

Case 3.1.3: $8n-14\leq s_0\leq s_1\leq 8n-7$.

In this case, we have $2(8n-14)\leq s_0+s_1\leq |S|\leq 12n-13$. It implies that $n=3$. Thus, $|S|\leq 12n-13=23$ and $10\leq s_0\leq s_1\leq 17$. Also, we have $s_c\leq|S|-s_0-s_1\leq 23-2\times 10=3$. Note that every vertex in $AQ_{n,3}^0$ (resp.\ $AQ_{n,3}^1$) has two distinct extra neighbors in $AQ_{n,3}^2$, at most one vertex in $(AQ_{n,3}^0-S_0)\cup (AQ_{n,3}^1-S_1)$ can be disconnected with $AQ_{n,3}^2-S_2$.
If a vertex $v$ in $AQ_{n,3}^0$ (resp.\ $AQ_{n,3}^1$) remains a singleton in $AQ_{n,3}-F$, then all the extra edges incident with $v$ are in $S_c$. It implies that $s_c\geq 4$, a contradiction. Thus, any component of $AQ_{n,3}^0-S_0$ (resp.\ $AQ_{n,3}^1-S_1$) is connected to $AQ_{n,3}^2-S_2$.
It leads to that $AQ_{n,3}-S$ is connected, and so $|V(H)|=|V(AQ_{n,3})|$.

Case 3.2: $k\geq4$.

Without loss of generality, assume that $0\in I$. We consider the following two cases according to the value of $I\setminus\{0\}$.

Case 3.2.1: $I=\{0,1\}$ or $I=\{0,k-1\}$.

Without loss of generality, assume that $I=\{0,1\}$ and $s_0\leq s_1$. For $j\in [k]_0\setminus\{0,1\}$, $AQ_{n,k}^j-S_j$ is connected. By Lemma~\ref{aq1}(3), since $s_c\leq|S|-s_0-s_1\leq (12n-13)-2(4n-6)=4n-1<2k^{n-1}$ for $n\geq 3$ and $k\geq 4$, there exists a fault-free edge joining $AQ_{n,k}^j-S_j$ and $AQ_{n,k}^{j+1}-S_{j+1}$ for $j\in [k]_0\setminus\{0,1,k-1\}$. It leads to $AQ_{n,k}^J-S_J$ is connected. By Lemma~\ref{aq1}(2), since $|S|\leq 12n-13$, we have $4n-6\leq s_0\leq s_1\leq |S|-(4n-6)\leq 8n-7$. The following three cases should be considered.

Case 3.2.1a: $4n-6\leq s_0\leq s_1\leq 8n-15$.

Case 3.2.1b: $4n-6\leq s_0\leq 8n-15$ and $8n-14\leq s_1\leq8n-7$.

Case 3.2.1c: $8n-14\leq s_0\leq s_1\leq 8n-7$.

Since the discussions for Case 3.2.1a, Case 3.2.1b and Case 3.2.1c are similar to those of Case 3.1.1, Case 3.1.2 and Case 3.1.3, respectively, the details are omitted.

Case 3.2.2: $I=\{0,t\}$, where $2\leq t\leq k-2$.

For $1\leq m\leq t-1$ or $t+1\leq m\leq k-1$, $AQ_{n,k}^m-S_m$ is connected.
By the similar argument as Case~3.2.1, those $(AQ_{n,k}^m-S_m)$'s for $1\leq m\leq t-1$ and $t+1\leq m\leq k-1$ belong to the components, say $C_1$ and $C_2$, respectively, of $AQ_{n,k}-S$.
Without loss of generality, assume that $s_0\leq s_t$. We consider the following cases.

Case 3.2.2a: $4n-6\leq s_0\leq s_t\leq 8n-15=8(n-1)-7$.

In this case, $AQ_{n,k}^0-S_0$ (resp.\ $AQ_{n,k}^t-S_t$) is disconnected. By Lemma~\ref{edge1}, $AQ_{n,k}^0-S_0$ (resp.\ $AQ_{n,k}^t-S_t$) has exactly two components: a large component, say $H_0$ (resp.\ $H_t$), and a singleton. Note that every vertex in $AQ_{n,k}^0$ has exactly two distinct extra neighbors in $AQ_{n,k}^1$ (resp.\ $AQ_{n,k}^{k-1}$). Since $s_c\leq|S|-s_0-s_t\leq (12n-13)-2(4n-6)=4n-1<2k^{n-1}-4$ for $n\geq 3$ and $k\geq 4$,
there is at least a fault-free edge between $H_0$ and $AQ_{n,k}^1-S_1$ (resp.\ $AQ_{n,k}^{k-1}-S_{k-1}$). This implies that
$H_0$ is connected to $AQ_{n,k}^1-S_1$ (which is part of $C_1$) and is connected to $AQ_{n,k}^{k-1}-S_{k-1}$ (which is part of $C_2$). It follows that $H_0$ is contained in both $C_1$ and $C_2$. By a similar discussion, $H_t$ is contained in both $C_1$ and $C_2$. This implies that $AQ_{n,k}-S$ has a large component $H=C_1=C_2$ and smaller components which contain at most two vertices in total. It leads to $|V(H)|\geq |V(AQ_{n,k})|-2$.

Case 3.2.2b: $4n-6\leq s_0\leq 8n-15$ and $8n-14\leq s_t\leq8n-7$.

In this case, we have $s_c\leq|S|-s_1-s_t\leq (12n-13)-(4n-6)-(8n-14)=7$. Since every vertex in $AQ_{n,k}^0$ (resp.\ $AQ_{n,k}^t$) has four distinct extra neighbors in $AQ_{n,k}-AQ_{n,k}^0-AQ_{n,k}^t$, any component with more than two vertices in $AQ_{n,k}^0-S_0$ is connected to $AQ_{n,k}^1-S_1$ (which is part of $C_1$) and is connected to $AQ_{n,k}^{k-1}-S_{k-1}$ (which is part of $C_2$). This implies that only a component with a singleton in $AQ_{n,k}^0-S_0$ can be disconnected with both $C_1$ and $C_2$. By a similar discussion, only a component with a singleton in $AQ_{n,k}^t-S_t$ can be disconnected with both $C_1$ and $C_2$. This implies that $AQ_{n,k}-S$ has a large component $H=C_1=C_2$ and smaller components which contain at most two vertices in total. It leads to $|V(H)|\geq |V(AQ_{n,k})|-2$.

Case 3.2.2c: $8n-14\leq s_0\leq s_t\leq 8n-7$.

In this case, we have $2(8n-14)\leq s_0+s_t\leq |S|\leq 12n-13$. It implies that $n=3$. Thus, $|S|\leq 12n-13=23$ and $10\leq s_0\leq s_t\leq 17$. Also, we have $s_c\leq|S|-s_0-s_t\leq 23-2\times 10=3$. Since every vertex in $AQ_{n,k}^0$ has four distinct extra neighbors in $AQ_{n,k}-AQ_{n,k}^0-AQ_{n,k}^t$, any component of $AQ_{n,k}^0-S_0$ is connected to $AQ_{n,k}^1-S_1$ (which is part of $C_1$) and is connected to $AQ_{n,k}^{k-1}-S_{k-1}$ (which is part of $C_2$). By a similar discussion, any component of $AQ_{n,k}^t-S_t$ is connected to $AQ_{n,k}^{t-1}-S_{t-1}$ (which is part of $C_1$) and is connected to $AQ_{n,k}^{t+1}-S_{t+1}$ (which is part of $C_2$). This implies that $AQ_{n,k}-S$ is connected, and so $|V(H)|=|V(AQ_{n,k})|$.

Case 4: $|I|=3$.

By Lemma~\ref{aq1}(2), for each $i\in I$, $s_i\geq 4(n-1)-2=4n-6$ and $s_i\leq (12n-13)-2(4n-6)=4n-1<8n-7$.
Since $|S|\leq 12n-13$, we have $s_c\leq |S|-3(4n-6)=5$. For each $i\in I$, $AQ_{n,k}^i-S_i$ is disconnected, and by Lemma~\ref{edge1}, $AQ_{n,k}^i-S_i$ has two components, one is the large component, say $H_i$, and the other is a singleton, say $v_i$.
Let $M$ be the union of smaller components of $AQ_{n,k}-S$.
We consider the following two cases according to $k=3$ or not.

Case 4.1: $k=3$.

In this case, $I=\{0,1,2\}$ and $J=\emptyset$. Since $s_c\leq 5<2k^{n-1}-4$ for $n\geq 3$, all $H_i$'s for $i\in I$ belong to the same component (i.e., $H$) in $AQ_{n,3}-S$. Clearly, $V(M)\subseteq \{v_0,v_1,v_2\}$. We claim $|V(M)|\leq 2$. Otherwise, $s_c\geq |N_{AQ_{n,3}^1\cup AQ_{n,3}^2}(v_0)|+|N_{AQ_{n,3}^2}(v_1)|=4+2=6$, a contradiction. This implies that $AQ_{n,3}-S$ has a large component and smaller components which contain at most two vertices in total. It leads to $|V(H)|\geq |V(AQ_{n,k})|-2$.

Case 4.2: $k\geq 4$.

We consider the following three cases.

Case 4.2.1: The three integers of $I$ are consecutive.

Without loss of generality, assume that $I=\{0,1,2\}$. For $j\in [k]_0\setminus\{0,1,2\}$, $AQ_{n,k}^j-S_j$ is connected. By Lemma~\ref{aq1}(3), there are $2k^{n-1}$ edges between subgraphs $AQ_{n,k}^j$ and $AQ_{n,k}^{j+1}$ $j\in[k]_0$. Since $s_c\leq 5<2k^{n-1}-4$ for $n\geq 3$ and $k\geq 4$,
all $H_i$'s for $i\in I$ and all subgraphs ($AQ_{n,k}^j-S_j$)'s for $j\in J$ belong to the same component (i.e., $H$) in $AQ_{n,k}-S$. Clearly, $V(M)\subseteq \{v_0,v_1,v_2\}$. By the similar discussion as Case 4.1, we can show that $|V(M)|\leq 2$ and $|V(H)|\geq |V(AQ_{n,k})|-2$.

Case 4.2.2: Only two integers in $I$ are consecutive.

Without loss of generality, assume that $I=\{0,1,t\}$, where $t\in\{3,4,\ldots,k-2\}$.
For $2\leq m\leq t-1$ or $t+1\leq m\leq k-1$, $AQ_{n,k}^m-S_m$ is connected.
By the similar argument as Case 4.2.1,
those $(AQ_{n,k}^m-S_m)$'s for $2\leq m\leq t-1$ and $t+1\leq m\leq k-1$ belong to the components, say $C_1$ and $C_2$, respectively, of $AQ_{n,k}-S$.
Since $s_c\leq 5<2k^{n-1}-4$ for $n\geq 3$ and $k\geq 4$,
$H_0$ is connected to $H_1$ and $AQ_{n,k}^{k-1}-S_{k-1}$ (which is part of $C_2$),
$H_1$ is connected to $AQ_{n,k}^2-S_2$ (which is part of $C_1$),
$H_t$ is connected to $AQ_{n,k}^{t-1}-S_{t-1}$ (which is part of $C_1$) and is connected to $AQ_{n,k}^{t+1}-S_{t+1}$ (which is part of $C_2$).
This implies that $AQ_{n,k}-S$ has a large component $H=C_1=C_2$ and $V(M)\subseteq \{v_0,v_1,v_t\}$.
We claim $|V(M)|\leq 2$. Otherwise, $s_c\geq |N_{AQ_{n,k}^{k-1}}(v_0)|+|N_{AQ_{n,k}^{t-1}\cup AQ_{n,k}^{t+1}}(v_t)|=2+4=6$, a contradiction.
This implies that $AQ_{n,k}-S$ has a large component $H$ and smaller components which contain at most two vertices in total. It leads to $|V(H)|\geq |V(AQ_{n,k})|-2$.

Case 4.2.3: None any two integers of $I$ are consecutive.

Without loss of generality, suppose $I=\{0,t,p\}$, where $2\leq t<p\leq k-2$ and $p-t\geq 2$.
For $1\leq m\leq t-1$ or $t+1\leq m\leq p-1$ or $p+1\leq m\leq k-1$, $AQ_{n,k}^m-S_m$ is connected.
By the similar argument as Case 4.2.1, those $(AQ_{n,k}^m-S_m)$'s for $1\leq m\leq t-1$, $t+1\leq m\leq p-1$, and $p+1\leq m\leq k-1$ belong to the components, say $C_1$, $C_2$ and $C_3$, respectively, of $AQ_{n,k}-S$. Since $s_c\leq 5<2k^{n-1}-4$ for $n\geq 3$ and $k\geq 4$,
$H_0$ is connected to $AQ_{n,k}^1-S_1$ (which is part of $C_1$) and is connected to $AQ_{n,k}^{k-1}-S_{k-1}$ (which is part of $C_3$),
$H_t$ is connected to $AQ_{n,k}^{t-1}-S_{t-1}$ (which is part of $C_1$) and is connected to $AQ_{n,k}^{t+1}-S_{t+1}$ (which is part of $C_2$),
$H_p$ is connected to $AQ_{n,k}^{p-1}-S_{p-1}$ (which is part of $C_2$) and is connected to $AQ_{n,k}^{p+1}-S_{p+1}$ (which is part of $C_3$).
This implies that $AQ_{n,k}-S$ has a large component $H=C_1=C_2=C_3$ and $V(M)\subseteq \{v_0,v_t,v_p\}$.
Since $4|V(M)|\leq s_c\leq 5$, we have $|V(M)|\leq 1$.
This implies that $AQ_{n,k}-S$ has a large component $H$ and smaller component which contain at most one vertex. It leads to $|V(H)|\geq |V(AQ_{n,k})|-1$.
\qed

\begin{theorem}\label{th2}
Let $AQ_{n,k}$ be the augmented $k$-ary $n$-cube,
where $n\geq 2$ and $k\geq 3$ are integers. Then
$AQ_{n,k}$ is $(4n-4)$-strongly Menger edge connected.
\end{theorem}
\pf
Let $F$ be an arbitrary faulty edge set of $AQ_{n,k}$ with
$|F|\leq 4n-4$. Since $\lambda(AQ_{n,k})=4n-2$, $AQ_{n,k}-F$ is connected.
Let $u,v\in V(AQ_{n,k})$ be any two distinct vertices such that $\text{deg}_{AQ_{n,k}-F}(u)\leq\text{deg}_{AQ_{n,k}-F}(v)$, and let $d_u=\text{deg}_{AQ_{n,k}-F}(u)$. From Proposition~\ref{m}, we need to show that the minimum size of a $(u,v)$-edge cut is equal to $d_u$. That is, we will show that $u$ and $v$ are still connected after the removal of at most $d_u-1$ edges in $AQ_{n,k}-F$.

Suppose, on the contrary, that $u$ and $v$ are separated by deleting a set of edges $E_f$ with $|E_f|\leq d_u-1$ in $AQ_{n,k}-F$. That is, $u$ is disconnected with $v$ in $AQ_{n,k}-(F\cup E_f)$. Since $d_u={\rm deg}_{AQ_{n,k}-F}(u)\leq {\rm deg}_{AQ_{n,k}}(u)=4n-2$, we have $|E_f|\leq 4n-3$. Let $S=F\cup E_f$. Then $|S|\leq (4n-4)+(4n-3)=8n-7$. By Lemma~\ref{edge1}, there is a component $H$
in $AQ_{n,k}-S$ such that $|V(H)|\geq |V(AQ_{n,k})|-1$. Clearly, $|V(H)|\ne|V(AQ_{n,k})|$, for otherwise, both $u$ and $v$ are contained in $H$.
Since $u$ is disconnected from $v$ in $AQ_{n,k}-S$, without loss of generality, assume that $u$ is a singleton in $AQ_{n,k}-S$.
Clearly, $E(\{u\},N_{AQ_{n,k}-F}(u))\subseteq E_f$. Thus, $|E_f|\geq|N_{AQ_{n,k}-F}(u)|=\text{deg}_{AQ_{n,k}-F}(u)=d_u$, which contradict to $|E_f|\leq d_u-1$. This shows that $AQ_{n,k}$ is $(4n-4)$-strongly Menger edge connected.
\qed

\begin{figure}[htb]
\begin{center}
\includegraphics[width=3.6in]{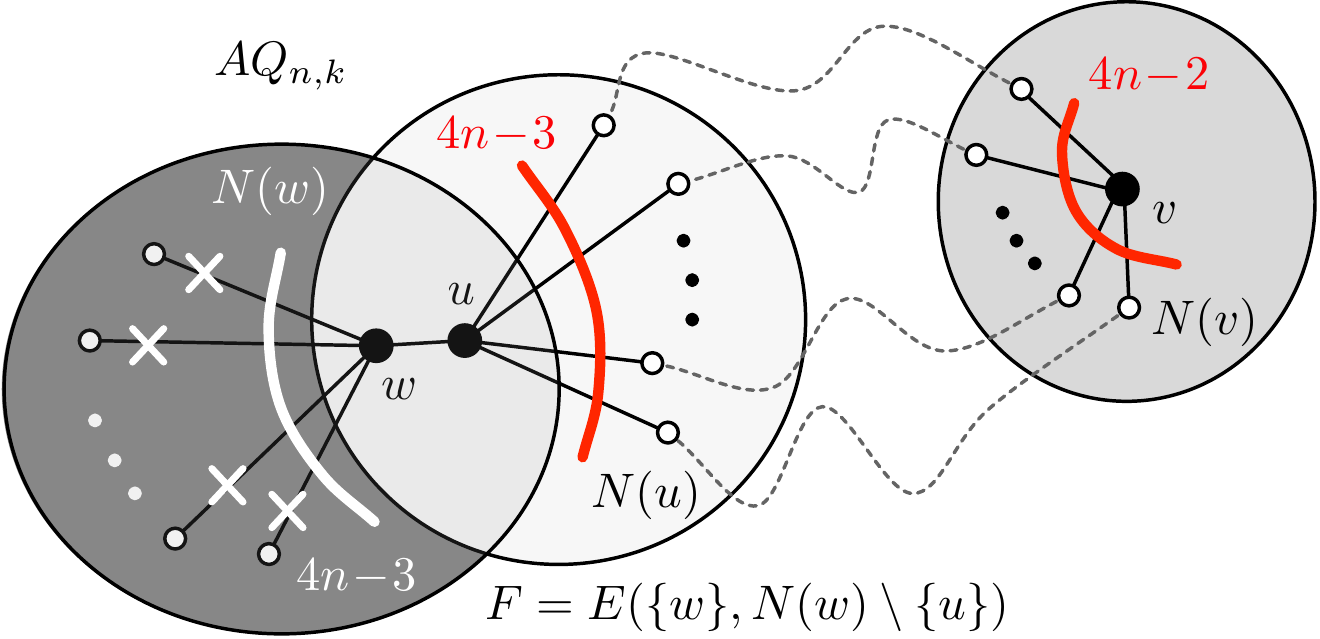}
\caption{Illustration for Remark~\ref{rem3}}
\label{remark3}
\end{center}
\end{figure}

\begin{remark}\label{rem3}
{\rm
To show that $AQ_{n,k}$ is not $(4n-3)$-strongly Menger edge connected, we consider the following example. See Fig.~\ref{remark3} for an illustration. Let $(u,w)\in E(AQ_{n,k})$ and $v\in V(AQ_{n,k})\setminus N[w]$. Let $F=E(\{w\},N(w)\setminus\{u\})$ be a faulty subset of edges in $AQ_{n,k}$ (i.e., edges with cross marks in the darkest area of Fig.~\ref{remark3}). Clearly, $|F|=4n-3$ and there are no more than $4n-3$ edge-disjoint paths between $u$ and $v$ in $AQ_{n,k}-F$. Since $\text{deg}_{AQ_{n,k}-F}(u)=\text{deg}_{AQ_{n,k}-F}(v)=4n-2$, $AQ_{n,k}$ is not $(4n-3)$-strongly Menger edge connected. Thus, the result of Theorem~\ref{th2} is optimal in the sense that the number of faulty edges cannot be increased.
}
\end{remark}

\begin{theorem}\label{th3}
Let $AQ_{n,k}$ be the augmented $k$-ary $n$-cubes,
where $n\geq 2$ and $k\geq 3$ are integers. Then
$AQ_{n,k}$ is $(8n-10)$-conditional strongly Menger edge connected.
\end{theorem}
\pf
Let $F$ be an arbitrary conditional faulty edge set of $AQ_{n,k}$ with $|F|\leq 8n-10$. Then $\delta(AQ_{n,k}-F)\geq 2$. Let $u$ and $v$ be any two distinct vertices in $AQ_{n,k}-F$ such that $\text{deg}_{AQ_{n,k}-F}(u)\leq\text{deg}_{AQ_{n,k}-F}(v)$. Also, let $d_u=\text{deg}_{AQ_{n,k}-F}(u)$ and let $d_v=\text{deg}_{AQ_{n,k}-F}(v)$. From Proposition~\ref{m}, we will show that $u$ and $v$ are connected by $d_u$ edge-disjoint fault-free paths in $AQ_{n,k}-F$. This means that $u$ and $v$ are still connected if the number of edges deleted is no more than $d_u-1$ in $AQ_{n,k}-F$.

Suppose, on the contrary, that $u$ and $v$ are separated by deleting a set of edges $E_f$ with $|E_f|\leq d_u-1\leq d_v-1$ in $AQ_{n,k}-F$. Let $S=F\cup E_f$. That is, $u$ and $v$ are disconnected in $AQ_{n,k}-S$. Since $d_u={\rm deg}_{AQ_{n,k}-F}(u)\leq {\rm deg}_{AQ_{n,k}}(u)=4n-2$, we have $|E_f|\leq 4n-3$.  Then $|S|\leq (8n-10)+(4n-3)=12n-13$. By Lemma~\ref{edge2}, there is a component $H$ in $AQ_{n,k}-S$ such that $|V(H)|\geq |V(AQ_{n,k})|-2$. It means that there are at most two vertices in $AQ_{n,k}-S$ not belonging to $H$. Clearly, $|V(H)|\ne|V(AQ_{n,k})|$, for otherwise, both $u$ and $v$ are contained in $H$.

If $|V(H)|=|V(AQ_{n,k})|-1$, without loss of generality, assume $u$ is a singleton in $AQ_{n,k}-S$. Clearly, $E(\{u\},N_{AQ_{n,k}-F}(u))\subseteq E_f$. Thus, $|E_f|\geq|N_{AQ_{n,k}-F}(u)|=\text{deg}_{AQ_{n,k}-F}(u)=d_u$, which contradict to $|E_f|\leq d_u-1$.

Assume that $|V(H)|=|V(AQ_{n,k})|-2$. Let $x$ and $y$ be the two vertices which are not belonging to $H$ in $AQ_{n,k}-S$. Consider the following two cases:

Case 1: $x$ and $y$ are adjacent in $AQ_{n,k}-S$.

Since $u$ and $v$ are separated in $AQ_{n,k}-S$, without loss of generality, we assume that $u\in V(H)$ and $v=x$. Clearly, $E(\{x,y\},N_{AQ_{n,k}-F}(\{x,y\}))\subseteq E_f$. Since $F$ is a conditional faulty edge set, $d_{AQ_{n,k}-F}(y)\geq 2$, thus there is at least one edge except $(x,y)$ which is incident with $y$ in $E_f$. Thus, $|E_f|\geq |N_{AQ_{n,k}-F}(\{x,y\})|\geq(\text{deg}_{AQ_{n,k}-F}(x)-1)+1=d_v$, which contradicts to $|E_f|\leq d_v-1$.

Case 2: $x$ and $y$ are not adjacent in $AQ_{n,k}-S$.

Since $u$ and $v$ are separated in $AQ_{n,k}-S$, we have $\{u,v\}\cap \{x,y\}\neq \emptyset$. Without loss of generality, assume $u=x$. Clearly, $E(\{u\},N_{AQ_{n,k}-F}(u))\subseteq E_f$. Thus, $|E_f|\geq|N_{AQ_{n,k}-F}(u)|=\text{deg}_{AQ_{n,k}-F}(u)=d_u$, which contradicts to $|E_f|\leq d_u-1$.
\qed

\begin{figure}[htb]
\begin{center}
\includegraphics[width=3.6in]{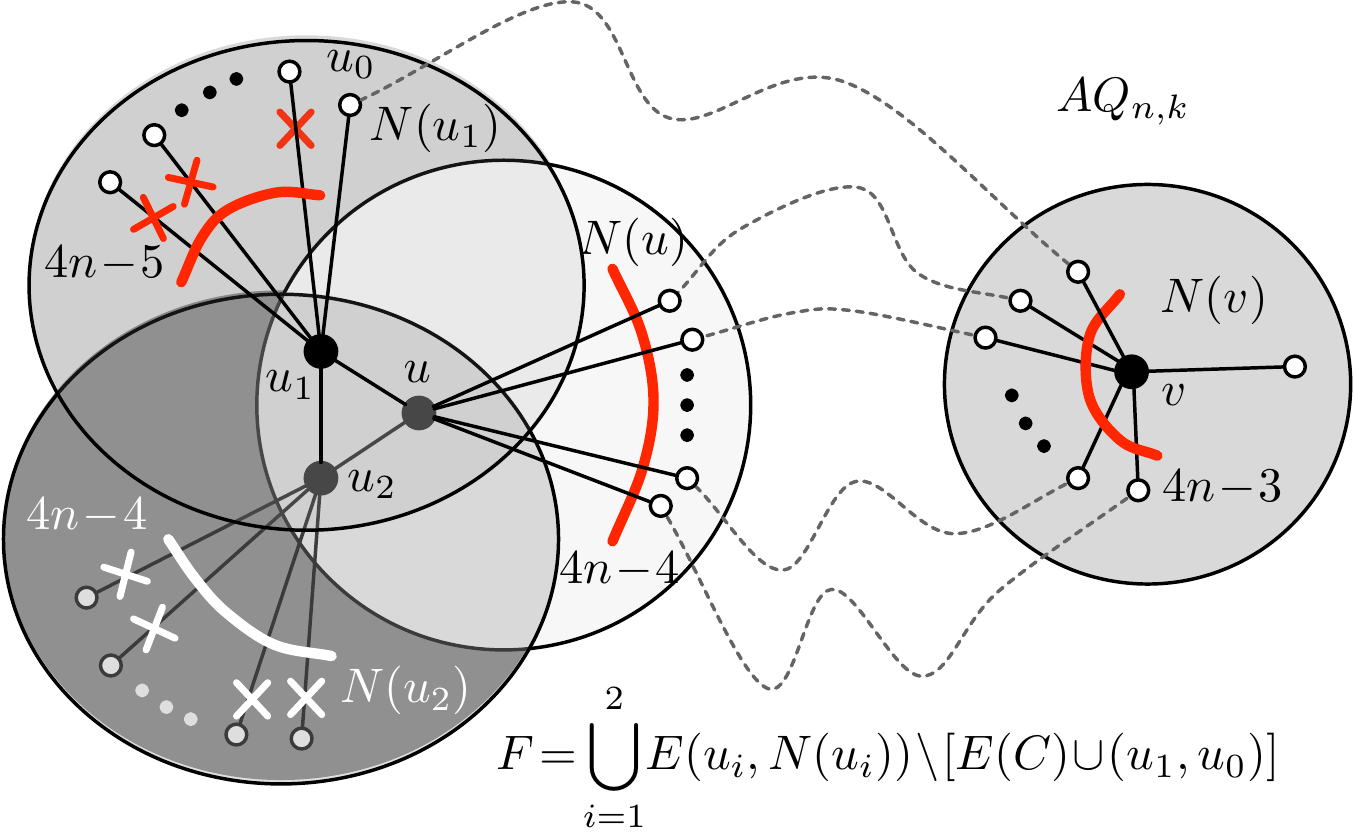}
\caption{Illustration for Remark~\ref{rem4}}
\label{remark4}
\end{center}
\end{figure}

\begin{remark}\label{rem4}
{\rm
To show that $AQ_{n,k}$ is not $(8n-9)$-conditional strongly Menger edge connected, we consider the following example. See Fig.~\ref{remark4} for an illustration. Let $C=(u,u_1,u_2,u)$ be a $3$-cycle in $AQ_{n,k}$. Also, let $u_0\in N(u_1)\setminus\{u,u_2\}$ and $v\in V(AQ_{n,k})\setminus(N(u_1)\cup N(u_2))$. Let $F=\bigcup_{i=1}^2 E(u_i,N(u_i))\setminus[E(C)\cup(u_1,u_0)]$ be a faulty subset of edges in $AQ_{n,k}$ (i.e., edges with cross marks in Fig.~\ref{remark4}). Clearly, $|F|=2(4n-2)-1-4=8n-9$ and there are no more than $4n-3$ edge-disjoint paths between $u$ and $v$ in $AQ_{n,k}-F$. Since $\text{deg}_{AQ_{n,k}-F}(u)=\text{deg}_{AQ_{n,k}-F}(v)=4n-2$, $AQ_{n,k}$ is not $(8n-9)$-conditional strongly Menger edge connected. Thus, the result of Theorem~\ref{th3} is optimal in the sense that the number of faulty edges cannot be increased.
}
\end{remark}

\section{Concluding remarks}
In literature, there are many papers with results of computing strong Menger (edge) connectivity in several popular classes of triangle-free graphs. In this paper, we study the strong Menger (edge) connectivity of one kind of graph which has many triangles, namely augmented $k$-ary $n$-cube $AQ_{n,k}$. By exploring and utilizing the structural properties of $AQ_{n,k}$, we show that $AQ_{n,3}$ (resp.\ $AQ_{n,k}$, $k\geq 4$) is $(4n-9)$-strongly (resp.\ $(4n-8)$-strongly) Menger connected for $n\geq 3$, and $AQ_{n,k}$ is $(4n-4)$-strongly Menger edge connected for $n\geq 2$ and $k\geq 3$.
Moreover, under the restricted condition that each vertex has at least two fault-free edges, we show that $AQ_{n,k}$ is $(8n-10)$-conditional strongly Menger edge connected for $n\geq 2$ and $k\geq 3$.
All results we obtained are optimal in the sense of the maximum number of tolerated vertex (resp.\ edge) faults. Intuitively, we think that this method can also be applied to other $r$-regular $r$-connected graphs with triangles.

\section*{Appendix A}

{\it
Let $S$ be an arbitrary set of edges in $AQ_{2,k}$ for $k\geq 3$. If $|S|\leq 9$, then there exists a component $H$ in $AQ_{2,k}-S$ such that $|V(H)|\geq |V(AQ_{2,k})|-1$.
}

\pf Let $H$ be the large component of $AQ_{2,k}-S$. By Lemma~\ref{aq1}(2), the result holds if $|S|\leq \kappa(AQ_{2,k})-1=5$. Now assume $6\leq |S|\leq 9$.
Note that $AQ_{2,k}$ contains $k$ disjoint copies of $k$-cycle, say $AQ_{2,k}^i$, $i\in[k]_0$. Let $S_i=S\cap E(AQ_{2,k}^i)$ and $s_i=|S_i|$ for $i\in [k]_0$.
Let $I=\{i\in [k]_0\colon\,AQ_{2,k}^i-S_i\ \text{is disconnected}\}$ and $J=[k]_0\setminus I$.
Clearly, $s_i\geq 2$ for each $i\in I$. So $|I|\leq \min\{k,4\}$. Let $M$ be the union of smaller components of $AQ_{n,k}-S$. We consider the following cases.

Case 1. $|I|=0$.

For $j\in[k]_0$, $AQ_{2,k}^j-S_j$ is connected. By Lemma~\ref{aq1}(3), there are $2k$ edges between adjacent subgraphs $AQ_{2,k}^j$ and $AQ_{2,k}^{j+1}$ for $j\in[k]_0$. Since $k\geq 3$ and $|S|\leq 9<2\times(2k)$, there exists at most one integer, say $i\in[k]_0$, such that all the edges between $AQ_{2,k}^i$ and $AQ_{2,k}^{i+1}$ are faulty. Since there is a fault-free edge between $AQ_{2,k}^j-S_j$ and $AQ_{2,k}^{j+1}-S_{j+1}$ for each $j\in[k]_0\setminus\{i\}$, it implies that $AQ_{2,k}-S$ is connected. Thus, $|V(H)|=|V(AQ_{2,k})|$.

Case 2. $|I|=1$.

Without loss of generality, assume that $I=\{0\}$. So, $V(M)\subset V(AQ_{2,k}^0)$. By Lemma~\ref{aq1}(2), $s_0\geq 4(n-1)-2=2$. For $j\in J$, $AQ_{2,k}^j-S_j$ is connected.
By Lemma~\ref{aq1}(3), there are $2k$ edges between adjacent subgraphs $AQ_{2,k}^j$ and $AQ_{2,k}^{j+1}$ for $j\in[k]_0$. Since $s_c\leq|S|-s_0\leq 9-2=7<2k$ for $k\geq 4$. If $k\geq 4$, there exists a fault-free edge between $AQ_{2,k}^j-S_j$ and $AQ_{2,k}^{j+1}-S_{j+1}$ for all $j\in[k]_0$. It leads to $AQ_{2,k}^J-S_J$ is connected. Since every vertex in $AQ_{2,k}^0$ has four distinct extra neighbors, we have $4|V(M)|\leq s_c\leq 7$. Thus, $|V(M)|\leq 1$ and $|V(H)|\geq |V(AQ_{2,k})|-1$.

If $k=3$ (see Fig.~\ref{aq23-24-33}(a)), then $AQ_{2,3}^0-S_0$ is either three singletons or an edge together with a singleton. Note that every vertex in $AQ_{2,3}^0$ has four distinct extra neighbors and $2\leq s_0\leq 3$. Since $s_0+4|V(M)|\leq 9$, we have $|V(M)|\leq 1$. Thus, at least a vertex of $AQ_{2,3}^0-S_0$ is connected to $AQ_{2,3}^1-S_1$ and $AQ_{2,3}^2-S_2$. It implies $|V(H)|\geq |V(AQ_{2,3})|-1$.

Case 3. $|I|=2$.

Without loss of generality, assume that $I=\{0,t\}$ where $t\in[k]_0\setminus\{0\}$.
For $j\in J$, $AQ_{2,k}^j-S_j$ is connected.
For $i\in I$, it is clear that $s_i\geq 2$.
Let $H_i$ be the large component of $AQ_{2,k}^i-S_i$.
Since $s_c\leq |S|-s_0-s_t\leq 9-2\times 2=5$ and every vertex in $AQ_{2,k}^i$ has four distinct extra neighbors, $AQ_{2,k}^i-S_i$ has two components, one is $H_i$ and the other is a singleton, say $v_i$.
By Lemma~\ref{aq1}(3), there are $2k$ edges between adjacent subgraphs $AQ_{2,k}^j$ and $AQ_{2,k}^{j+1}$ for $j\in[k]_0$. Since $s_c\leq 5<2k$ for $k\geq 3$, all $H_i$'s for $i\in I$ and all subgraphs ($AQ_{2,k}^j-S_j$)'s for $j\in J$ belong to the same component (i.e., $H$) in $AQ_{2,k}-S$.
Since the two singletons, $v_0$ in $AQ_{2,k}^0-S_0-H_0$ and $v_t$ in $AQ_{2,k}^t-S_t-H_t$, may be adjacent in $AQ_{2,k}-S$, we have $s_0+s_t+(4|V(M)|-2)\leq 9$. Thus, $|V(M)|\leq 1$ and $|V(H)|\geq |V(AQ_{2,k})|-1$.

Case 4. $3\leq |I|\leq \min\{k,4\}$.

In this case, $k\geq 3$.
For $j\in J$, $AQ_{2,k}^j-S_j$ is connected.
For $i\in I$, it is clear that $s_i\geq 2$.
Let $H_i$ be the large component of $AQ_{2,k}^i-S_i$.
By Lemma~\ref{aq1}(3), there are $2k$ edges between adjacent subgraphs $AQ_{2,k}^j$ and $AQ_{2,k}^{j+1}$ for $j\in[k]_0$. Since $s_c\leq |S|-2|I|\leq 9-2\times 3=3<2k$, all $H_i$'s for $i\in I$ and all subgraphs ($AQ_{2,k}^j-S_j$)'s for $j\in J$ belong to the same component (i.e., $H$) in $AQ_{2,k}-S$.
Since $s_c\leq 3$ and every vertex in $AQ_{2,k}^i$ has four distinct extra neighbors, if a vertex of $AQ_{2,k}^i-S_i$ is a singleton, then it must be connected to $H$. Thus, $|V(H)|\geq |V(AQ_{2,k})|$.
\qed

\section*{Appendix B}

{\it
Let $S$ be an arbitrary set of edges in $AQ_{2,k}$ for $k\geq 3$. If $|S|\leq 11$, then there exists a component $H$ in $AQ_{2,k}-S$ such that $|V(H)|\geq |V(AQ_{2,k})|-2$.
}

\pf Let $n=2$ and $H$ be the large component of $AQ_{2,k}-S$. By Lemma~\ref{edge1}, the result holds if $|S|\leq 8n-7=9$. Now we consider $10\leq |S|\leq 11$.
Recall that $AQ_{2,k}$ contains $k$ disjoint copies of $k$-cycle, say $AQ_{2,k}^i$, $i\in[k]_0$. Let $S_i=S\cap E(AQ_{2,k}^i)$ and $s_i=|S_i|$ for $i\in [k]_0$.
Let $I=\{i\in [k]_0\colon\,AQ_{2,k}^i-S_i\ \text{is disconnected}\}$ and $J=[k]_0\setminus I$.
Clearly, $s_i\geq 2$ for each $i\in I$. So $|I|\leq \min\{k,5\}$. Let $M$ be the union of smaller components of $AQ_{n,k}-S$. We consider the following cases.

Case 1. $|I|=0$.

For $j\in[k]_0$, $AQ_{2,k}^j-S_j$ is connected. By Lemma~\ref{aq1}(3), there are $2k$ edges between adjacent subgraphs $AQ_{2,k}^j$ and $AQ_{2,k}^{j+1}$ for $j\in[k]_0$. Since $k\geq 3$ and $|S|\leq 11<2\times(2k)$, there exists at most one integer, say $i\in[k]_0$, such that all the edges between $AQ_{2,k}^i$ and $AQ_{2,k}^{i+1}$ are faulty. Since there is a fault-free edge between $AQ_{2,k}^j-S_j$ and $AQ_{2,k}^{j+1}-S_{j+1}$ for each $j\in[k]_0\setminus\{i\}$, it implies that $AQ_{2,k}-S$ is connected. Thus, $|V(H)|=|V(AQ_{2,k})|$.

Case 2. $|I|=1$.

Without loss of generality, assume that $I=\{0\}$. So, $V(M)\subset V(AQ_{2,k}^0)$. By Lemma~\ref{aq1}(2), $s_0\geq 4(n-1)-2=2$. For $j\in J$, $AQ_{2,k}^j-S_j$ is connected.
By Lemma~\ref{aq1}(3), there are $2k$ edges between adjacent subgraphs $AQ_{2,k}^j$ and $AQ_{2,k}^{j+1}$ for $j\in[k]_0$. Since $s_c\leq|S|-s_0\leq 11-2=9<2k$ for $k\geq 5$. If $k\geq 5$, there exists a fault-free edge between $AQ_{2,k}^j-S_j$ and $AQ_{2,k}^{j+1}-S_{j+1}$ for all $j\in[k]_0$. It leads to $AQ_{2,k}^J-S_J$ is connected. Since every vertex in $AQ_{2,k}^0$ has four distinct extra neighbors, we have $4|V(M)|\leq s_c\leq 9$. Thus, $|V(M)|\leq 2$ and $|V(H)|\geq |V(AQ_{2,k})|-2$.

If $k=3$ (see Fig.~\ref{aq23-24-33}(a)), then $AQ_{2,3}^0-S_0$ is either three singletons or an edge together with a singleton. Note that every vertex in $AQ_{2,3}^0$ has four distinct extra neighbors and $2\leq s_0\leq 3$. Since $s_0+4|V(M)|\leq 11$, we have $|V(M)|\leq 2$. Thus, at least a vertex of $AQ_{2,3}^0-S_0$ is connected to $AQ_{2,3}^1-S_1$ and $AQ_{2,3}^2-S_2$. It implies $|V(H)|\geq |V(AQ_{2,3})|-2$.

If $k=4$ (see Fig.~\ref{aq23-24-33}(b)), then $AQ_{2,4}^0-S_0$ is one of the following: (1) four singletons; (2) an edge and two singletons; (3) a $2$-path (i.e., a path of length 2) and a singleton; (4) two nonadjacent edges. Note that every vertex in $AQ_{2,4}^0$ has four distinct extra neighbors and $2\leq s_0\leq 4$. Since $s_0+4|V(M)|\leq 11$, we have $|V(M)|\leq 2$, and $M$ is either an edge or at most two singletons. Thus, there exists either an edge, a 2-path, or two singletons of $AQ_{2,4}^0-S_0$, which and all subgraphs ($AQ_{2,4}^j-S_j$)'s for $j\in J$ belong to the same component (i.e., $H$) in $AQ_{2,k}-S$.
This shows that $|V(H)|\geq |V(AQ_{2,4})|-2$.

%
%

Case 3. $|I|=2$.

Without loss of generality, assume that $I=\{0,t\}$ where $t\in[k]_0\setminus\{0\}$.
For $j\in J$, $AQ_{2,k}^j-S_j$ is connected.
For $i\in I$, it is clear that $s_i\geq 2$.
Let $H_i$ be the large component of $AQ_{2,k}^i-S_i$.
Since $s_c\leq |S|-s_0-s_t\leq 11-2\times 2=7$ and every vertex in $AQ_{2,k}^i$ has four distinct extra neighbors, $AQ_{2,k}^i-S_i$ has two components, one is $H_i$ and the other is a singleton, say $v_i$.
By Lemma~\ref{aq1}(3), there are $2k$ edges between adjacent subgraphs $AQ_{2,k}^j$ and $AQ_{2,k}^{j+1}$ for $j\in[k]_0$. Since $s_c\leq 7<2\times(2k)$ for $k\geq 3$, all $H_i$'s for $i\in I$ and all subgraphs ($AQ_{2,k}^j-S_j$)'s for $j\in J$ belong to the same component (i.e., $H$) in $AQ_{2,k}-S$.
Since the two singletons, $v_0$ in $AQ_{2,k}^0-S_0-H_0$ and $v_t$ in $AQ_{2,k}^t-S_t-H_t$, may be adjacent in $AQ_{2,k}-S$, we have $s_0+s_t+(4|V(M)|-2)\leq 11$. Thus, $|V(M)|\leq 2$ and $|V(H)|\geq |V(AQ_{2,k})|-2$.


Case 4. $|I|=3$.

Without loss of generality, assume that $I=\{0,t,p\}$ where $t,p\in[k]_0\setminus\{0\}$. For $j\in J$, $AQ_{2,k}^j-S_j$ is connected.
For $i\in I$, it is clear that $s_i\geq 2$.
Let $H_i$ be the large component of $AQ_{2,k}^i-S_i$.
Since $s_c\leq |S|-s_0-s_t-s_p\leq 11-3\times 2=5$ and every vertex in $AQ_{2,k}^i$ has four distinct extra neighbors, $AQ_{2,k}^i-S_i$ has two components, one is $H_i$ and the other is a singleton, say $v_i$.
By Lemma~\ref{aq1}(3), there are $2k$ edges between adjacent subgraphs $AQ_{2,k}^j$ and $AQ_{2,k}^{j+1}$ for $j\in[k]_0$. Since $s_c\leq 5<2\times(2k)$ for $k\geq 3$, all $H_i$'s for $i\in I$ and all subgraphs ($AQ_{2,k}^j-S_j$)'s for $j\in J$ belong to the same component (i.e., $H$) in $AQ_{2,k}-S$. If every two singletons of $\{v_0,v_t,v_p\}$ are adjacent in $AQ_{2,k}-S$, then $s_0+s_t+s_p+(4\times 3-6)\geq 12$, a contradiction. Thus at least two singletons in $\{v_0,v_t,v_p\}$ are nonadjacent, and it follows that $s_0+s_t+s_p+(4|V(M)|-4)\leq 11$. This shows that $|V(M)|\leq 2$ and $|V(H)|\geq |V(AQ_{2,4})|-2$.


Case 5. $4\leq |I|\leq \min\{k,5\}$.

In this case, $k\geq 4$. For $j\in J$, $AQ_{2,k}^j-S_j$ is connected.
For $i\in I$, it is clear that $s_i\geq 2$.
Let $H_i$ be the large component of $AQ_{2,k}^i-S_i$.
By Lemma~\ref{aq1}(3), there are $2k$ edges between adjacent subgraphs $AQ_{2,k}^j$ and $AQ_{2,k}^{j+1}$ for $j\in[k]_0$. Since $s_c\leq |S|-2|I|\leq 11-2\times 4=3<2k$, all $H_i$'s for $i\in I$ and all subgraphs ($AQ_{2,k}^j-S_j$)'s for $j\in J$ belong to the same component (i.e., $H$) in $AQ_{2,k}-S$.
Since $s_c\leq 3$ and every vertex in $AQ_{2,k}^i$ has four distinct extra neighbors, if a vertex of $AQ_{2,k}^i-S_i$ is a singleton, then it must be connected to $H$. Thus, $|V(H)|\geq |V(AQ_{2,k})|$.
\qed

%
%


\begin{thebibliography}{99} \small
\baselineskip=12pt
\parskip 0pt



\bibitem{cll15}
H.-Y. Cai, H.-Q. Liu and M. Lu,
Fault-tolerant maximal local-connectivity on Bubble-sort star graphs,
Discrete Appl. Math. 181 (2015) 33--40.

\bibitem{cct14}
Y.-C. Chen, M.-H. Chen and J.J.M. Tan,
Maximally local connectivity and connected components of augmented cubes,
Inform. Sci. 273 (2014) 387--392.

\bibitem{cxu18} Q. Cheng, P. S. Li and M. Xu, Conditional (edge) fault-tolerant strong Menger (edge) connectivity of folded hypercubes,
Theoret. Comput. Sci. 728 (2018) 1--8.

\bibitem{cs02}
S.A. Choudum and V. Sunitha,
Augmented cubes,
Networks 40(2) (2002) 71--84.

\bibitem{gu16}
M.-M.~Gu, R.-X.~Hao, Y.-Q.~Feng,
The pessimistic diagnosability of bubble-sort star graphs and augmented $k$-ary $n$-cubes,
Int. J. Comput. Math.: Comput. Sys. Theory 1 (2016) 98--112.

\bibitem{hhc18} S. He, R.-X. Hao and E. Cheng,
Strongly Menger-edge-connectedness and strongly Menger-vertex-connectedness of regular networks,
Theoret. Comput. Sci. 731 (2018) 50--67.

\bibitem{lx18} P. Li and M. Xu,
Fault-tolerant strong Menger (edge) connectivity and 3-extra edge-connectivity
of balanced hypercubes, Theoret. Comput. Sci. 707 (2018) 56--68.

\bibitem{lx19}
P. Li and M. Xu,
Edge-fault-tolerant strong Menger edge connectivity on the class of hypercube-like networks,
Discrete Appl. Math. 259 (2019) 145--152.

\bibitem{lx19t}
P. Li and M. Xu, The $t/k$-diagnosability and strong Menger connectivity on star graphs with conditional faults, Theoret. Comput. Sci. 793 (2019) 181--192.


\bibitem{lz15}
R. Lin and H. Zhang,
The restricted edge-connectivity and restricted connectivity of augmented $k$-ary $n$-cubes,
Inter. J. Comput. Math. 93(8) (2016) 1281--1298.

\bibitem{Menger}
K. Menger,
Zur allgemeinen kurventheorie,
Fund. Math. 10 (1927) 95--115.

\bibitem{oc03}
E. Oh and J. Chen,
On strong Menger-connectivity of star graphs,
Discrete Appl. Math. 129 (2003) 499--511.

\bibitem{oh04}
E. Oh,
On strong fault tolerance (or strong Menger-connectivity) of multicomputer networks,
(Ph.D. thesis), Computer Science A\&M  University, August 2004.

\bibitem{qyang17}
Y. Qiao and W. Yang,
Edge disjoint paths in hypercubes and folded hypercubes with conditional faults,
Appl. Math. Comput. 294 (2017) 96--101.

\bibitem{sm19}
E. Sabir and J. Meng,
Parallel routing in regular networks with faults,
Inform. Proces. lett. 142 (2019) 84--89.

\bibitem{lscc}
L. Shih, C. Chiang, L. Hsu and J. Tan,
Strong Menger connectivity with conditional faults on the class of hypercube-like networks,
Inform. Process. Lett. 106 (2008) 64--69.

\bibitem{scht09}
C.-M. Shih, C.-F. Chiang, L.-H. Hsu and J. Tan,
Fault-tolerant maximal local-connectivity on Cayley graphs generated by transposition trees,
J. Interconnect. Netw. 10 (2009) 253--260.

\bibitem{wz}
S. Wang and N. Zhao,
The two-good-neighbor connectivity
and diagnosability of the augmented
three-ary $n$-cubes, The Computer journal, doi:10.1093/comjnl/bxy125


\bibitem{xs11}
Y. Xiang and I. A. Stewart,
Augmented $k$-ary $n$-cubes,
Inform. Sci. 181 (2011) 239--256.

\bibitem{xs09}
Y. Xiang and I. A. Stewart,
Pancyclicity and panconnectivity in augmented $k$-ary $n$-cubes,
IEEE 15th International Conference on Parallel and Distributed Systems (ICPADS),
Shenzhen, China, 2009, pp. 308--315.

\bibitem{yang17}
W. Yang, S. Zhao and S. Zhang,
Strong Menger connectivity with conditional faults of folded hypercubes,
Inform. Process. Lett. 125 (2017) 30--34.

\end{thebibliography}
\end{document}